\documentclass[final,12pt]{elsarticle}
\usepackage{subcaption}
\usepackage{amsmath}
\usepackage{amsbsy} 
\usepackage{tabularx}
\usepackage{stmaryrd}
\usepackage{enumitem}
\usepackage{siunitx}
\usepackage[numbers]{natbib}
\usepackage{amssymb}
\usepackage{xcolor}
\usepackage{geometry}
\makeatletter
\def\ps@pprintTitle{%
  \let\@oddhead\@empty
  \let\@evenhead\@empty
  \def\@oddfoot{} % Clear the odd footer completely
  \def\@evenfoot{} % Clear the even footer completely
}
\makeatother

\begin{document}

\begin{frontmatter}



\title{Physically Recurrent Neural Networks for Computational Homogenization of Composite Materials with Microscale Debonding}


\author[a,b,c]{N. Kovács}
\author[a]{M.A. Maia}
\author[a]{I.B.C.M. Rocha}
\author[b,c]{C. Furtado}
\author[b,c]{P.P. Camanho}
\author[a]{F.P. van der Meer}

\address[a]{{Delft University of Technology, Department of Civil Engineering and Geosciences}, {P.O.Box 5048}, 2600GA, Delft, {The Netherlands}}
\address[b]{{DEMec, Faculdade de Engenharia, Universidade do Porto}, {Rua Dr. Roberto Frias}, 4200-465 Porto, {Portugal}}
\address[c]{{INEGI, Instituto de Ciência e Inovação em Engenharia Mecânica e Engenharia Industrial}, {Rua Dr. Roberto Frias, 400, 4200-465 Porto}, {Portugal}}

\begin{abstract}
The growing use of composite materials in engineering applications has accelerated the demand for computational methods to accurately predict their complex behavior. Multiscale modeling based on computational homogenization is a potentially powerful approach for this purpose, but its widespread adoption is prevented by its excessive computational costs. A popular approach to address this computational bottleneck is using surrogate models, which have been used to successfully predict a wide range of constitutive behaviors. However, applications involving microscale damage and fracture remain largely unexplored. This work aims to extend a recent surrogate modeling approach, the Physically Recurrent Neural Network (PRNN), to include the effect of debonding at the fiber-matrix interface while capturing path-dependent behavior. The core idea of the PRNN is to implement the exact material models from the micromodel into one of the layers of the network. In this work, additional material points with a cohesive zone model are integrated within the network, along with the bulk points associated to the fibers and/or matrix. The limitations of the existing architecture are discussed and taken into account for the development of novel architectures that better represent the stress homogenization procedure. In the proposed layout, the history variables of cohesive points act as extra latent features that help determine the local strains of bulk points. Different architectures are evaluated starting with small training datasets. To maximize the predictive accuracy and extrapolation capabilities of the network, various configurations of bulk and cohesive points are explored, along with different training dataset types and sizes.
\end{abstract}


% \begin{highlights}
% \end{highlights}


\end{frontmatter}


\section{Introduction}
\label{}
Composite materials such as fiber-reinforced polymers (FRPs) can be tailored to reach optimal mechanical properties and enhance structural performance by design at different length scales. In practice, however, while modeling techniques have the potential to reduce the experimental testing costs and accelerate product development, the multi-scale nature and nonlinear behavior of these materials pose challenges for the reliable prediction of their global response. A suitable method to tackle such challenges is computational homogenization, or FE$^2$, a concurrent multiscale approach that offers a high-fidelity representation of the material behavior. Nevertheless, this flexibility and generality is followed by excessively high computational costs.

Several methods have been proposed to alleviate the computational bottleneck associated with $\text{FE}^2$. One alternative is to solve the equilibrium problem at the microscale, where a Representative Volume Element (RVE) is used to describe the material geometry and properties, in low-dimensional spaces for instance with Proper Orthogonal Decomposition (POD) \cite{MONTEIRO2008704} or Proper Generalized Decomposition (PGD) \cite{ELHALABI2013183,CREMONESI2013275}. Although both techniques can effectively reduce computational time, obtaining a meaningful set of basis functions/modes is not trivial. For instance, the construction of the basis functions based on the separation of the variables (\emph{e.g.}, space, time, geometry, etc.) in PGD requires careful formulation and involves a large number of parameters as non-linearity increases. On the other hand, POD relies on a meaningful coverage of potential loading paths in the snapshots used for mode extraction, which can become a difficult task when dealing with path-dependent materials. % TODO MARINA: mor: modes depend a lot on the snapshot selection, difficult when path dependent materials are considered cause even more snapshots will be needed (cite iuris paper for that)
% pgd: lots of parameters, more general than pod but , relies on the variables being separable, not as widespread as pod (simpler) 

% TODO MARINA: change the intro of this paragraph to specify what problem, and link the RNNS to something with less physical meaning/grounds, but more general
With higher potential to accelerate the microscopic problem, another popular approach is to replace the full-order microscale BVP with a \emph{surrogate model}. In \cite{GHAVAMIAN2019112594} and \cite{LOGARZO2021113482}, the full-order micromodel is replaced by a recurrent neural network (RNN) at each integration point. The network is then used to predict the average stress from the micromodel, as well as history variables such as the equivalent plastic strain. A drawback of this method is related to the black-box nature of these RNNs. Since the network only sees the input strain and output stress values during training with no physical interpretation, it is not able to extrapolate to unseen loading scenarios. This means that a large amount of data is required to train the network. 

A different approach was presented in a recent work using mechanics-informed machine learning \cite{LI2023117473}. In this approach, \textit{a priori} knowledge of mechanics was implemented in the network to predict elastoplasticity for composites with pressure-independent yielding behavior by physics-based decomposition of stresses and strains. The network showed accurate predictions trained on a small dataset. A newly emerging method of combining the mechanical basis with the flexibility of neural networks is using physics-augmented neural networks \cite{rosenkranz2024viscoelastictyphysicsaugmentedneuralnetworks,klein2023parametrisedpolyconvexhyperelasticityphysicsaugmented,schommartz2024physicsaugmentedneuralnetworksconstitutive}. The essence of these networks is that they are tailored for a specific material model to fulfil mechanical conditions by construction. Although these networks provide accurate predictions, their complexity increases significantly with more complex material models.

Another approach using physically recurrent neural networks (PRNNs) \cite{MAIA2023115934} is not limited to a specific type of material model. In PRNNs, the constitutive relations in the full-order micromodel are directly implemented in the hidden layer of the network, bypassing the need for a complex network while giving physical interpretation to the hidden layer. The performance of the PRNN was first assessed in several loading cases at the microscale and compared to state-of-the-art RNNs. Then, in a multiscale setting, the proposed network was used to replace the microsctructure with elastic and elastoplastic phases in the online evaluations \cite{MAIA2023115934}. Trained on monotonic data, the PRNN was able to predict with the same accuracy as the RNN, but with much less training data. Notably, the PRNN was able to extrapolate to non-monotonic test data even though it only saw data from monotonic loading scenarios during training. One of the trained PRNNs was applied in a $\text{FE}^2$ problem, and it was observed that replacing the full-order micromodel with the network resulted in a reduction of computational time by more than 20,000 times. %These exceptional results show great promise for further use of the PRNN and for extending it to different material models.

% Although a large part of the literature has been devoted to modeling the behavior of heterogeneous materials with path-dependency in the surrogate modeling field \cite{GHAVAMIAN2019112594,LOGARZO2021113482}, the task is still a quite challenging endeavor and a number of issues related to the black-box nature of conventional models are an issue. This is even more complicated once damage and fracture need to be considered in the micromodel.

A large part of recent literature on surrogate modeling focuses on predicting (hyper)elastic or elastoplastic behavior. Surrogate modeling for damage and fracture mechanics applications is a much less explored field. For example, Wang et al. \cite{WANG2019216} used deep reinforcement learning techniques to create traction-separation laws, but did not apply this in a multiscale setting. In the works of Liu \cite{LIU2020112913,LIU2021113914}, a deep material network (DMN) was developed, which describes the RVE with a network that is built up from physics-based building blocks. In the first work \cite{LIU2020112913}, debonding effects in the RVE were successfully captured by the adaptable cohesive building blocks included in the network. In a multi-stage training strategy, one of the phases from the DMN trained for the bulk material is enriched with the interface interaction learned from a second DMN built on top of cohesive building blocks. As a result, the accuracy of the final network is somehow limited by the original bulk DMN. In \cite{LIU2021113914} this method was extended to localization problems with a cell-division scheme, which overcame the difficulties related to selecting the proper size of the RVE. While these networks excel in extrapolating from linear elastic RVE data to nonlinear and path-dependent behavior, training and online evaluation are not straightforward. These two stages involve different input spaces and an iterative Newton-Raphson scheme is required for the prediction stage \cite{LIU2020112913,LIU2021113914}. Further, a probabilistic machine learning approach using Bayesian regression was proposed in \cite{ROCHA2021100083} and also applied to active learning of traction-separation relations in a multiscale setting, but this approach was not made suitable for capturing unloading behavior.

Even without a multiscale framework, the computational cost involved in fracture mechanics problems with FE can also be extreme. A noteworthy approach to deal with that is proposed by Kerfrieden et al. \cite{KERFRIDEN2013169} through the use of a domain separation strategy. In this work, this method is introduced to focus the computational power on the fracture, which requires most of the attention. The domain separation strategy can also used in multiscale settings, as proposed by Oliver et al. \cite{OLIVER2017560}. To further reduce the sampling points, it is combined with another key technique based on model order reduction, specifically the Reduced Optimal Quadrature (ROQ) . However, these methods are highly dependent on snapshots and the complexity of the problem increases with nonlinear and path-dependent materials.
% TODO Marina: link it to hyper reduction and mention same problems as before with pod based models: they depend a lot on snapshots and it becomes complicated with nonlinear and pathdependent materials

In this work, the ability of the PRNN to describe the effect of microscale damage is investigated. The study is restricted to diffuse damage in the form of microscale debonding in composite materials. The aim is to describe the stiffness degradation that is the result of diffuse damage, without the occurrence of macroscopic softening. Section 2 outlines the $\text{FE}^2$ framework and briefly describes the PRNN as detailed in \cite{MAIA2023115934}. Section 3 provides the data generation method for training and testing the networks. In section 4, the performance of the existing PRNN on cases with stiffness degradation is tested, motivating the development of new architectures. Section 5 introduces the proposed architectures, while Section 6 evaluates their performance in achieving the research objective.
% \newpage
\section{Theoretical background}
\label{}
In the following sections, the foundational aspects of the methods used in this work are discussed. This includes an overview of the FE$^2$ method and the homogenization procedure, as well as the main features and limitations of the existing PRNN.

\subsection{The $\text{FE}^\text{2}$ method}\label{fenn}

Computational homogenization with the $\text{FE}^2$ method allows for capturing the response of composite materials, for cases where it is difficult to do that on the macroscale due to the complex interaction between nonlinear constituents and microstructural geometry. In the $\text{FE}^2$ approach, the structure is discretized as a homogeneous macrostructure, where a heterogeneous micromodel is nested into each macroscopic integration point of it \cite{Schröder2014,FEYEL2000309,FEYEL20033233}. The micromodel is assumed to be a representative volume element (RVE). The macroscopic strain values are downscaled as boundary conditions for the micromodel, where the microscopic boundary value problem (BVP) is solved. The microscopic stress values obtained from the BVP are then upscaled back to the macromodel after a homogenization operation. This bypasses the need for any assumptions that would need to be made for formulating a macroscale constitutive relation.

\begin{figure}[!h]
\centering
\includegraphics[scale=0.64]{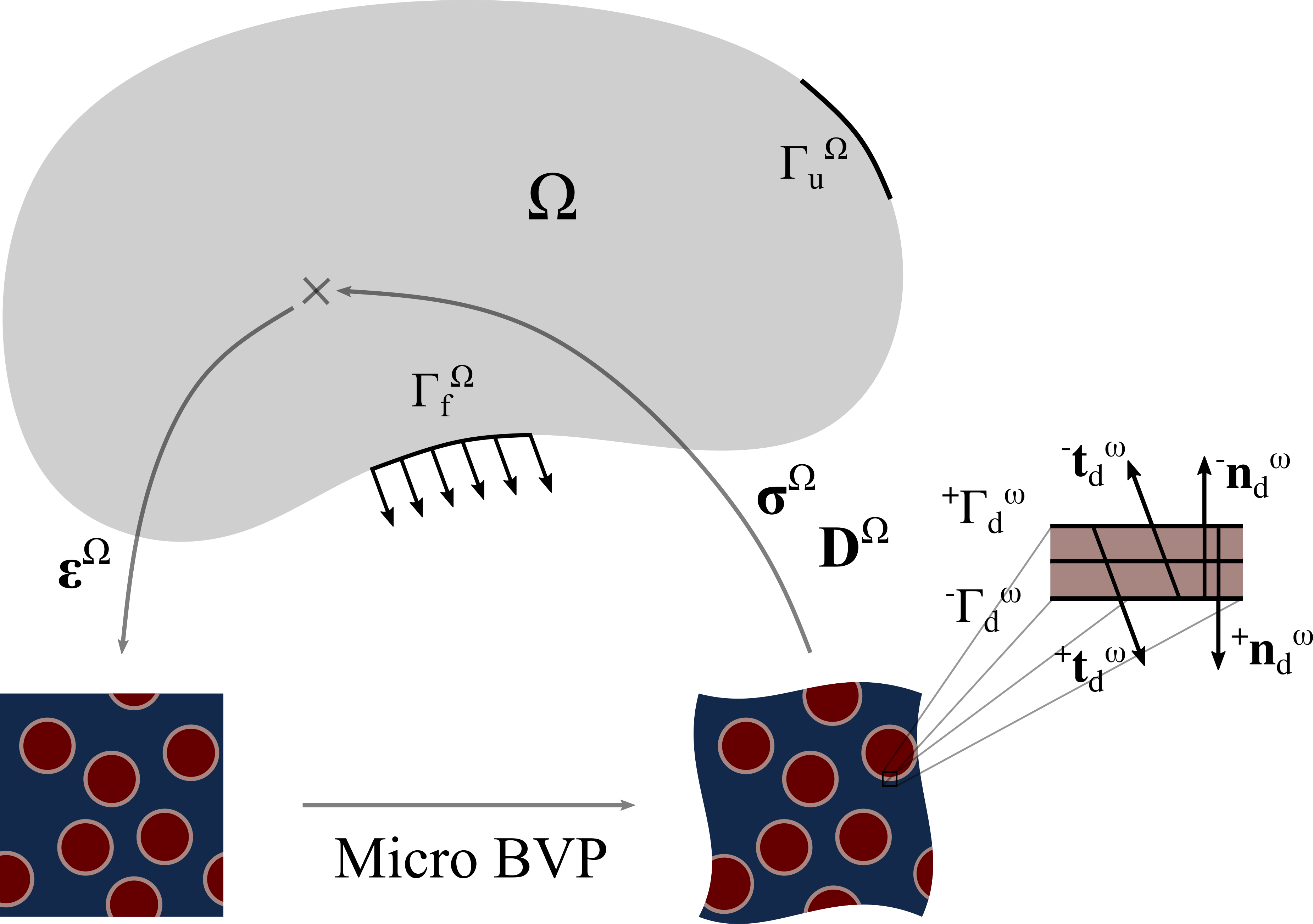}
\caption{\centering $\text{FE}^2$ framework}
\label{fig:fe2}
\end{figure}

The schematics of the $\text{FE}^2$ method is shown in Figure \ref{fig:fe2}. The macroscopic solid domain is denoted by $\Omega$, and the surfaces where the Dirichlet and Neumann boundary conditions are applied are denoted as $\Gamma_u^\Omega$ and $\Gamma_f^\Omega$, respectively. In this work, we simulate damage at the fiber-matrix interface, which precludes global softening of the micromodel and avoids the need for inserting a discontinuity on the macroscale. The discontinuity in the microscopic domain is denoted by $\Gamma_d^\omega$. At the fracture surface, the two opposite sides of the crack are differentiated by a $+$ and a $-$ sign. %The following key elements of $\text{FE}^2$ applied to damage mechanics are discussed in this chapter: the macroscale problem formulation, the microscale RVE problem formulation along with the boundary conditions received from the macroscale, and the micro- to macroscale transition. 

Instead of using a macroscale constitutive model, the macroscale strains at the integration points $\pmb{\varepsilon}^{\Omega}$ are passed to the RVE, where they are used as boundary conditions. The displacement field in the micromodel with domain denoted as $\omega$ is approximated for given boundary conditions with a finite element discretization of the RVE geometry. The macroscale strain is considered to be constant over the volume aside from a periodic microscale fluctuation field due to the assumption of separation of scales. Nonlinearity in the microscale problem comes from the constitutive models $\mathcal{D}^{\omega}$ and $\mathcal{T}^\omega$:
\begin{align}
\pmb{\sigma}^{\omega}, \pmb{\alpha}&=\mathcal{D}^{\omega}\left(\pmb{\varepsilon}^{\omega}, \pmb{\alpha}^{t-1} \right) \\
\mathbf{t}_d^\omega, d &= \mathcal{T}^\omega \left( \llbracket \mathbf{u}^\omega \rrbracket, d^{t-1}  \right),
\end{align}
where $\pmb{\sigma}^\omega$ is the microscale stress obtained from the microscale strain $\pmb{\varepsilon}^{\omega}$ and the internal variables $\pmb{\alpha}$, and $\mathbf{t}_d^\omega$ is the cohesive traction computed from the displacement jump $\llbracket \mathbf{u}^\omega \rrbracket$ and internal variable $d$.
After the computation of the full-order solution at the microscale, the resulting stress field is homogenized and the homogenized stress is returned to the macroscale model. The homogenized stress is defined as the average over the volume of the RVE:
\begin{equation}\label{eq:homogenization}
    \pmb{\sigma}^\Omega = \frac{1}{||\omega||} \int_\omega \pmb{\sigma}^\omega\, d\omega
\end{equation}

For accurate coupling between the two scales, the energy between them must be consistent. For this, the Hill-Mandel principle has to be satisfied \cite{HILL1963357} which is ensured when periodic boundary conditions are used.

To solve the macroscale problem stress update from $\pmb{\varepsilon}^{\omega}$ to $\pmb{\sigma}^{\omega}$ must be embedded in a nonlinear finite element solution procedure, e.g. based on the Newton-Raphson method. Note that the cost of solving the microscale BVP at each integration point and every macroscale iteration makes the FE$^2$ framework unfeasible for practical applications.

\subsection{Physically Recurrent Neural Network}

\begin{figure}[!h]
\minipage{0.49\textwidth}
  \centering
  \includegraphics[width=0.85\linewidth]{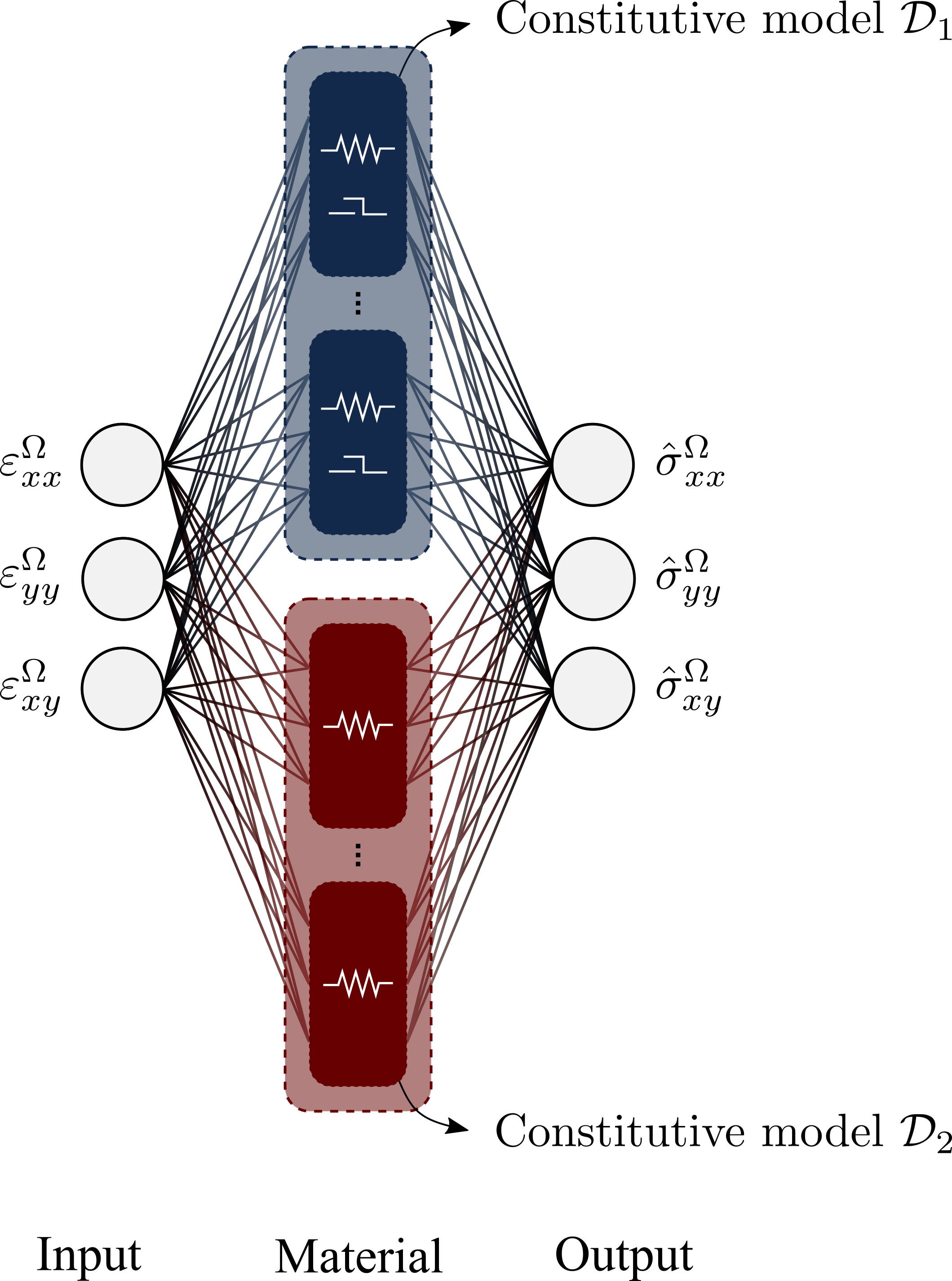}
\caption{\centering Physically recurrent neural network \cite{MAIA2023115934}} 
\label{fig:prnngen}
\endminipage\hfill
\minipage{0.49\textwidth}
  \centering
  \includegraphics[width=0.67\linewidth]{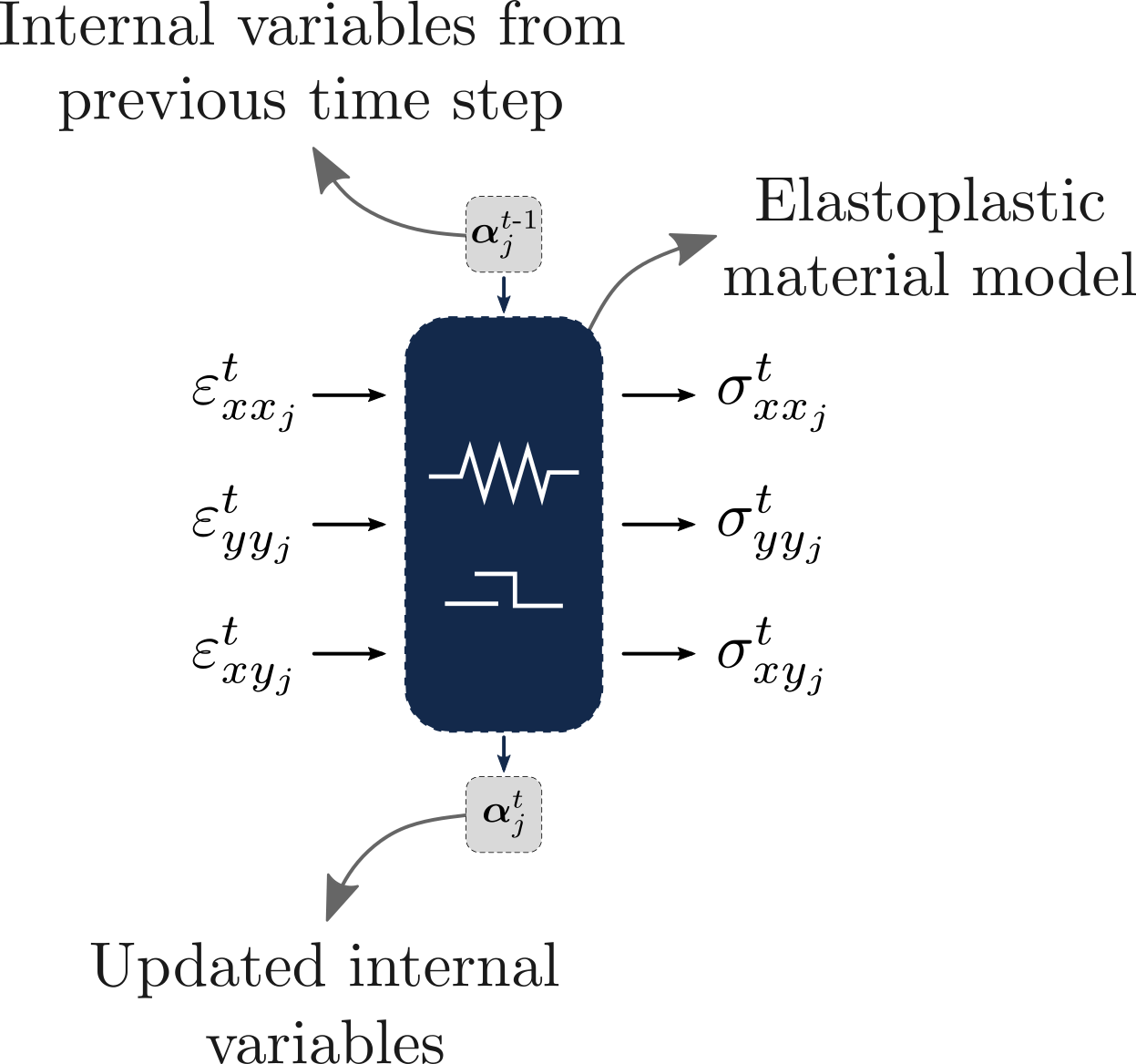}
\caption{\centering Neuron in material layer of physically recurrent neural network \cite{MAIA2023115934}}
\label{fig:prnncell}
\endminipage\hfill
\end{figure}

To tackle the issues related to the black-box nature of neural networks, Physically Recurrent Neural Networks (PRNNs) developed by Maia et al. \cite{MAIA2023115934} introduce a new way of embedding knowledge on the physics of a system in a surrogate for that system. Unlike in PINNs, where the physical constraints of the problem are incorporated in the loss function \cite{RAISSI2019686}, in PRNNs the actual material models used in the full-order microscopic BVP are implemented in the hidden layer of the network such that their state variables introduce a physics-based recurrency. Figure \ref{fig:prnngen} displays the PRNN in general terms for the case with two constitutive models on the microscale, $\mathcal{D}_1$ and $\mathcal{D}_2$.

The architecture consists of an input layer, a material layer, and an output layer. The macroscale strains $\pmb{\varepsilon}^\Omega$ at the integration points of the macrostructure are the inputs to this network. In two dimensions assuming small strains, this corresponds to 3 input values. These macroscale strains are passed through an encoder, which is a single dense layer with linear activation functions. This encoder converts the macroscale strains to microscale strains, which corresponds to the macro- to micro-scale transition in the $\text{FE}^2$ method. The microscale strains $\pmb{\varepsilon}$ are then given by:

\begin{equation}\label{eq:encoder}
    \pmb{\varepsilon} = \mathbf{W}_1 \pmb{\varepsilon}^\Omega + \mathbf{b}_1 ,
\end{equation}
where $\mathbf{W}_1$ are the weights connecting the input layer to the material layer and $\mathbf{b_1}$ is the bias associated with the encoder. There are no residual stresses considered in this work, which means that there is a zero stress state for when no strain is applied to the microstructure. Therefore, the network should also predict zero stresses when the strain inputs are zeros. This is achieved by setting the bias term $\mathbf{b_1} = \mathbf{0}$.

These microscale strains are passed through the material layer, which provides the essence of the physically recurrent neural network. The material layer consists of cells that each have three inputs and three outputs, and possibly internal variables. Each of these cells represents a fictitious integration point. In these points, a classical constitutive model $\mathcal{D}^\omega$ converts the microscale strains $\pmb{\varepsilon}$ to microscale stresses $\pmb{\sigma}$. This constitutive model is the exact same model that is used to compute the stress in the integration points of the full-order micromodel. The number of fictitious material points in the layer is as a model hyperparameter to be tweaked via model selection.

History dependence is naturally included in the PRNN by storing the internal variables $\pmb{\alpha}$ of each material point, which for example in plasticity can be plastic deformation, as they are computed in the assigned constitutive model. Therefore, path-dependency does not need to be learned from data. This stands in contrast with regular recurrent neural networks, where the evolution of history variables is also learned through additional learnable parameters and standard activation functions (e.g. with LSTM or GRU cells). The operation in the fictitious material point $j$ at time $t$ is shown in Figure \ref{fig:prnncell} and can be described by:

\begin{equation}\label{eq:matlay}
    \pmb{\sigma}_j^t, \pmb{\alpha}_j^t = \mathcal{D}^\omega(\pmb{\varepsilon}_j^t, \pmb{\alpha}_j^{t-1})
\end{equation}

After the local stresses are computed in the material layer, a decoder is applied to these values. This step is designed to mimic the stress homogenization step of the $\text{FE}^2$ method. Essentially, the local macroscale response is a combination of the response of the fictitious material points present in the material layer, with the network being tasked to learn how to best combine the response of a small number of material points into a representative macroscale response.

In the particular achitecture shown in Fig. \ref{fig:prnngen}, the decoder consists of a dense layer with a SoftPlus activation function applied on the weights. This is done to represent the homogenization process through numerical integration, in which weights are strictly positive. Therefore, the macroscale stress output of the network is obtained by:

\begin{equation}\label{eq:decoder}
    \widehat{\pmb{\sigma}}^\Omega = \phi_{\text{sp}}(\mathbf{W}_2) \pmb{\sigma} + \mathbf{b}_2 ,
\end{equation}
where $\mathbf{W}_2$ are the weights connecting the material layer to the output layer, and $\mathbf{b_2}$ is the bias associated with the decoder. This bias term is again set to zero to ensure zero macroscale stresses for zero local stress values.
% \begin{equation}
%     \widehat{\pmb{\sigma}}^\Omega = \phi_{sp}(\mathbf{W}_2) \pmb{\sigma}^\omega
% \end{equation}

During training, the following loss function is minimized:

\begin{equation}\label{eq:lossf}
    \text{L} = \frac{1}{N} \sum_{t=1}^{N} ||\pmb{\sigma}^\Omega (\pmb{\varepsilon}_t^\Omega) - \widehat{\pmb{\sigma}}^\Omega (\pmb{\varepsilon}_t^\Omega)||^2 ,
\end{equation}
where $N$ is the number of stress-strain pairs in the dataset and $\pmb{\sigma}^\Omega (\pmb{\varepsilon}_t^\Omega)$ is the target value. Stochastic gradient descent is used with Adam optimizer. Because of the history variables $\pmb{\alpha}$, back-propagation in time becomes necessary, stress/strain pairs are grouped in paths (time series) and gradients are only computed at the end of the paths. Data and gradient handling is therefore analogous to when training RNNs, with the key difference being that memory in PRNNs is physical and interpretable. It is worth emphasizing that if the material model is implemented with autodiff support, gradients are handled automatically through general-use packages such as pytorch and tensorflow. Otherwise, a detailed implementation of how to incorporate these using finite differences is given in \cite{MAIA2023115934}. 

In the full-order micromodel used in \cite{MAIA2023115934}, J$_2$ plasticity was used for the matrix and a linear elasticity model for the fibers. In the previous study, the PRNN was able to find a solution with only the elastoplastic model (i.e. the constitutive model used to describe the matrix) in the material layer. The expected linear elastic behavior in the fibers is reproduced in elastoplastic material points when small enough strain values are passed by the encoder and stresses are amplified in the decoder, making the matrix model effectively work like the linear elastic model of the fiber. 

% \begin{figure}[!h]
% \centering
% \includegraphics[scale=1.0]{Pictures/micro0.png}
% \caption{\centering Full-order micromodel used in \cite{MAIA2023115934}}
% \label{fig:micro0}
% \end{figure}
% This micromodel is subjected to different loading scenarios to obtain high-fidelity solutions used to train and assess the network's performance. 

\section{Data generation}
\label{}
To generate the data for assessing the PRNN's performance when predicting microscale damage, a micromodel with cohesive elements at the fiber-matrix interface is considered. This section introduces the micromodel used to create the training and test datasets and the different loading conditions that are considered.

\subsection{Full-order micromodel}\label{microm}

\begin{figure}[!h]
\centering
\includegraphics[scale=0.11]{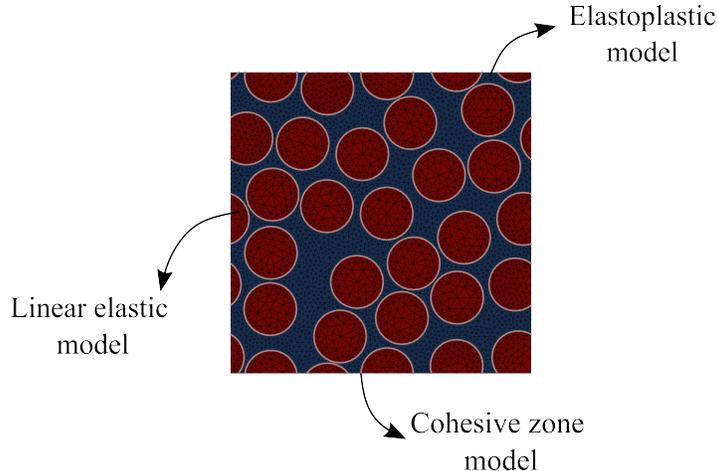}
\caption{\centering Full-order micromodel used in this work}
\label{fig:micro}
\end{figure}

The FE model of the microscale is shown in Figure \ref{fig:micro} and consists of 25 periodically arranged fibers with diameter of $\SI{5}{\micro\metre}$ embedded in a matrix to result in a fiber volume fraction of 0.6. There are two bulk constitutive models: a plasticity model for the matrix, and a linear elastic model for the fibers. The geometry, mesh, and bulk material properties in the RVE are kept as in \cite{MAIA2023115934} except that zero-thickness interface elements are positioned on the fiber-matrix interface. Limiting damage to the fiber-matrix interface means that no global failure can take place. Tractions in the interface elements are computed from the displacement jump with the bilinear cohesive zone model (CZM) by Turon et al. \cite{TURON20061072}.

% linear elastic fibers with properties $\text{E} = 74 000$ MPa and $\nu = 0.2$, and elastoplastic matrix modeled with J2 plasticity with properties $\text{E} = 3130$ MPa, $\nu = 0.3$
%     \sigma_y = 64.8 − 33.6 e^{-\varepsilon_{eq}^p/0.0003407}\\
%     \varepsilon_{eq}^p = \sqrt{\frac{2}{3} \pmb{\varepsilon}^p : \pmb{\varepsilon}^p}

For the CZM properties, we use equal normal and shear strength $\tau_n^0 = \tau_s^0 = 60$ MPa, mode I and mode II fracture energy $G_{Ic} = 0.874 \text{kJ/m}^2$, $G_{IIc} = 1.717 \text{kJ/m}^2$, mode interaction parameter $\eta = 1$, and penalty stiffness $K = 5 \cdot 10^7 \text{N/mm}^3$. Plane stress conditions are assumed for the micromodel.

\subsection{Load path generation}

To generate data for training and testing of the network, the micromodel is subjected to different loading paths using periodic boundary conditions. The datasets can be separated into two categories: proportional and non-proportional loading. 

\subsubsection{Proportional loading}

\begin{figure}[!h]
\begin{subfigure}{0.49\textwidth}
  \centering
  \includegraphics[width=1\linewidth]{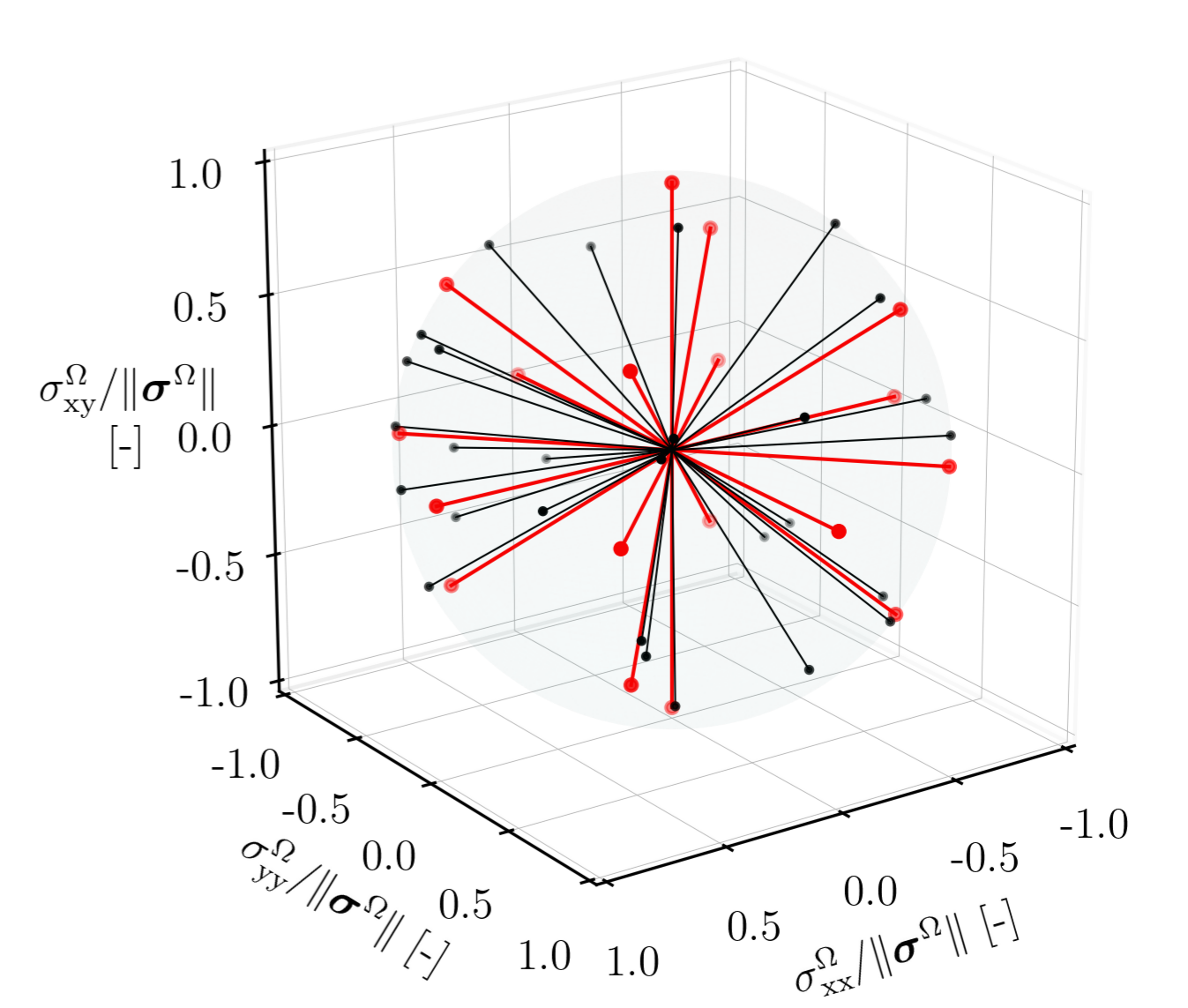}
\caption{\centering Proportional loading path directions}
\label{fig:cr}
\end{subfigure}
\hfill
\begin{subfigure}{0.49\textwidth}
\centering
\includegraphics[width=0.95\linewidth]{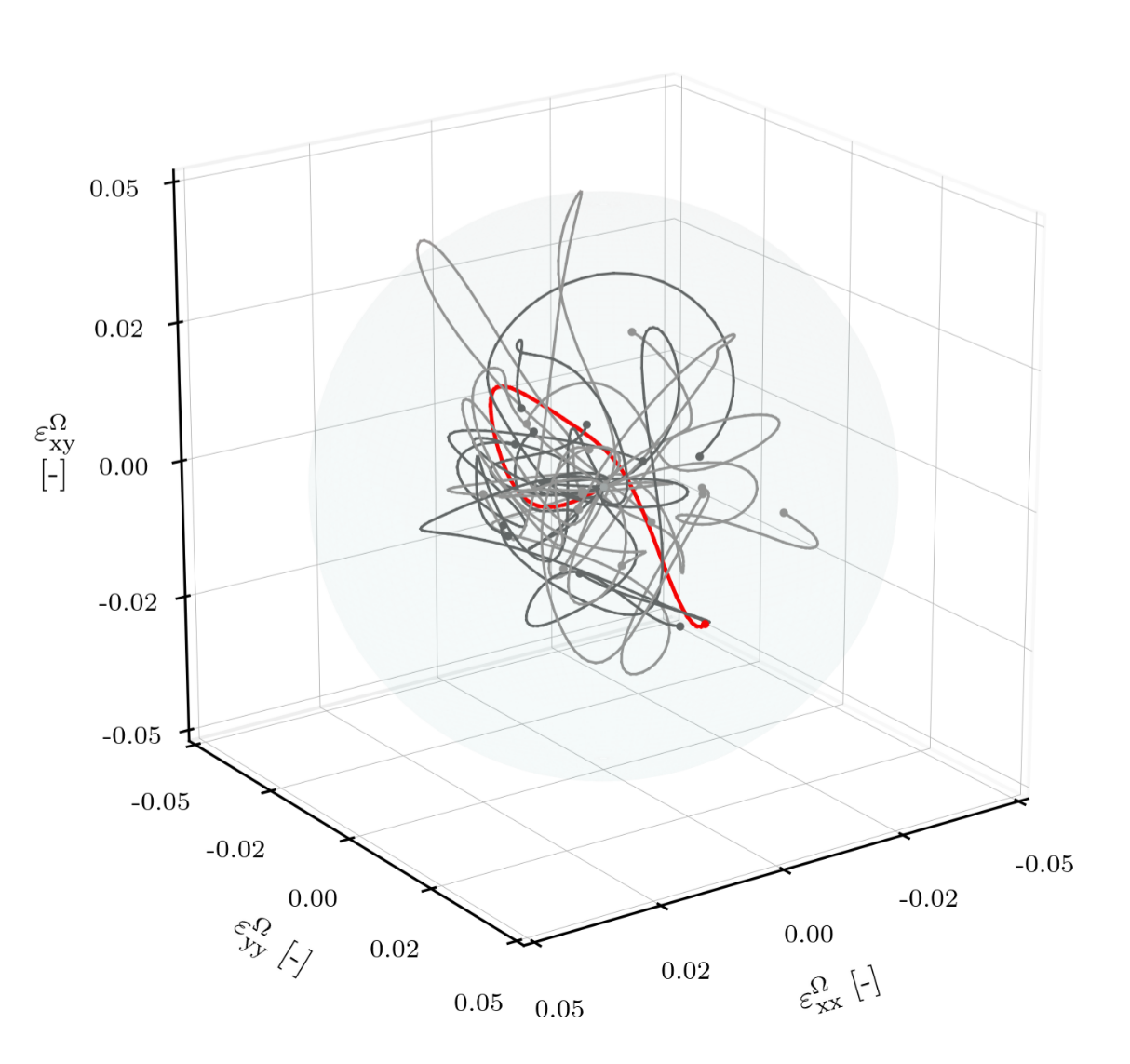}
\caption{\centering Non-proportional loading paths}
\label{fig:gp}
\end{subfigure}
\caption{Types of loading paths considered in this work.}
\end{figure}

For proportional loading, a modified arc-length algorithm is used to enforce the directions of the applied stress, as described in \cite{ROCHA2020103995}. In this method, the proportionality of the stress response is enforced by considering a constant load vector with increments defined in terms of displacement magnitudes, specifically the sum of the unsigned applied displacements. In this work, the increments are also fixed and worth $s = 1.67 \times 10^{-3}$ mm. The loading directions are categorized as either fundamental or random. The fundamental directions contain 18 common loading cases often used for material model calibration, shown in red in Figure \ref{fig:cr}, and include pure tension, compression, shear, biaxial tension, and combinations thereof. On the other hand, the random directions are obtained by sampling three values, each corresponding to a load vector component, from a normal distribution $\mathcal{N} (0, 1)$ and normalizing them to a unit vector, with examples shown in black in Figure \ref{fig:cr}.

\begin{figure}[h]
\centering
\begin{subfigure}{0.48\textwidth}
\centering
\includegraphics[width=\linewidth]{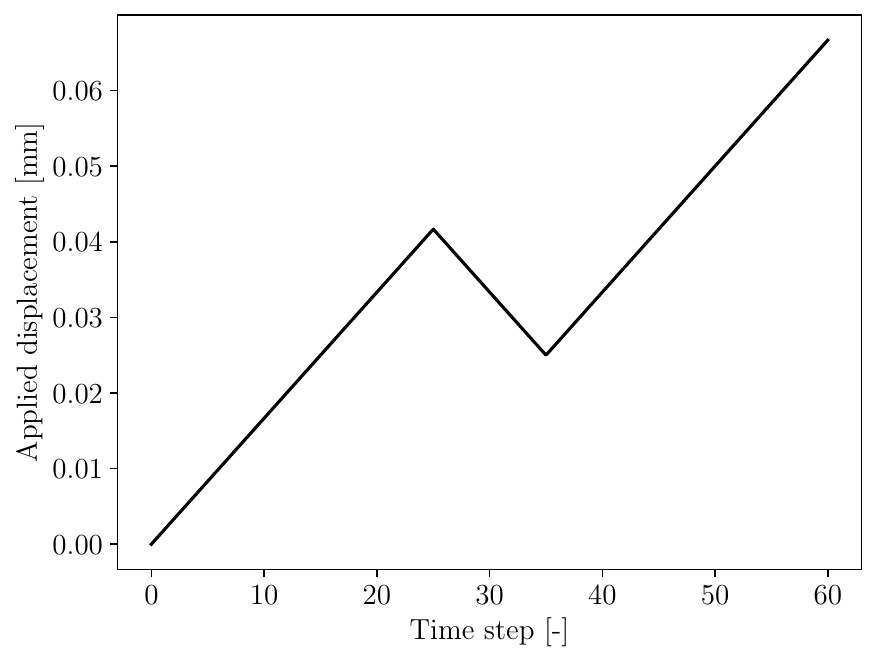}
\caption{\centering Function with one unloading cycle}\label{fig:load1}
\end{subfigure}
     \hfill
     \begin{subfigure}{0.48\textwidth}
\centering
  \includegraphics[width=\linewidth]{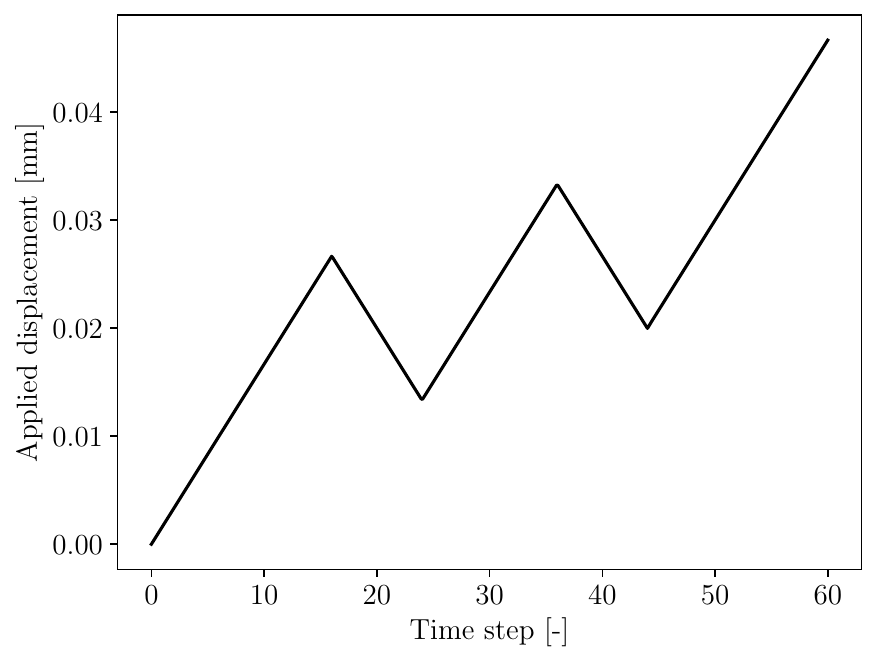}
  \caption{\centering Function with two unloading cycles}\label{fig:load2}
\end{subfigure}
        \caption{\centering Loading functions used to generate proportional loading curves}
        \label{fig:loadingc}
\end{figure}

% \begin{figure}[!h]
% \centering
% \includegraphics[scale=0.8]{Pictures/randomcurves1.png}
% \caption{\centering Proportional loading path directions}
% \label{fig:cr}
% \end{figure}

In this work, only non-monotonic loading is considered. During non-monotonic loading, while the direction in which the step size is kept fixed, unloading takes place at different loading steps for a predefined amount of time. The loading functions that define the relation between $t$ and the magnitude of loading for non-monotonic cases considered in this paper are shown in Figure \ref{fig:loadingc}. In the arc-length formulation, this corresponds to the imposed value for the unsigned sum of the displacements at the controlling nodes.

\subsection{Non-proportional loading}

\begin{figure}[!h]
    \centering
    \begin{subfigure}{0.49\textwidth}
\centering
\includegraphics[width=\linewidth]{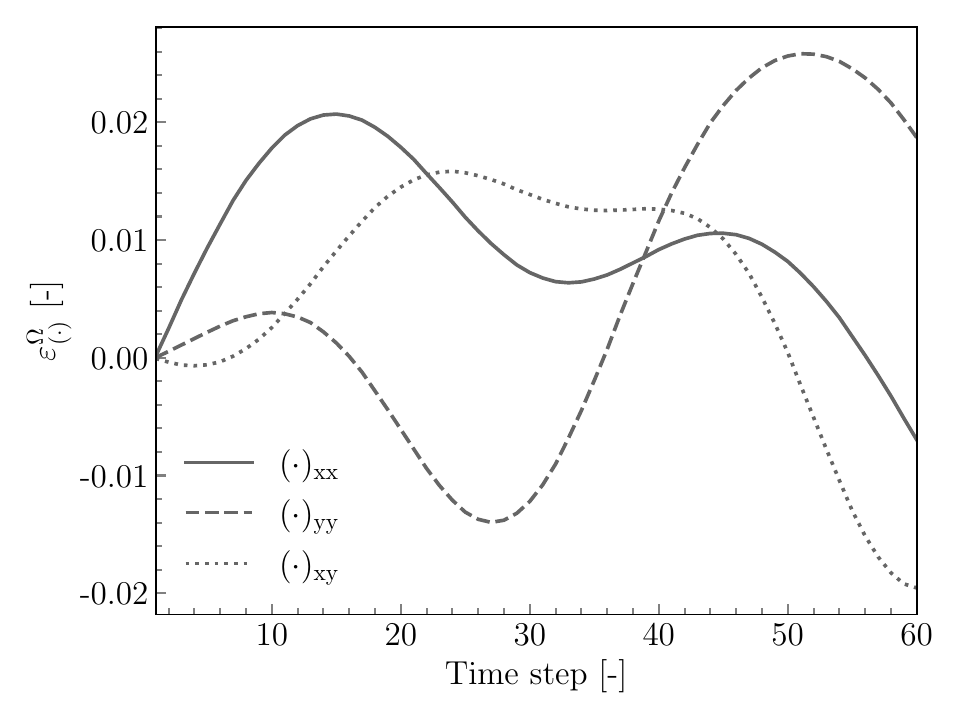}
\caption{\centering Strain paths}\label{fig:strains}
\end{subfigure}
     \hfill
     \begin{subfigure}{0.49\textwidth}
\centering
  \includegraphics[width=\linewidth]{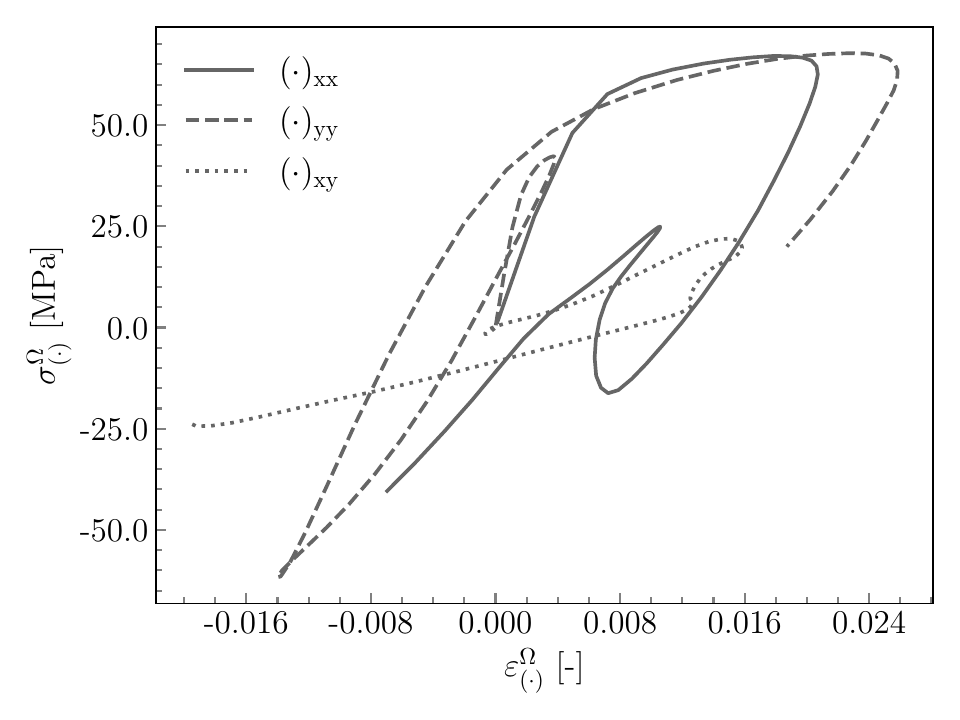}
  \caption{\centering Stress-strain curves}\label{fig:strstr}
\end{subfigure}
    \caption{\centering Example of non-proportional GP-based loading path}
    \label{fig:strains}
\end{figure}

To create more diverse loading scenarios, non-proportional and non-monotonic loading paths are generated. Both the direction of loading and the step size are varied at each time step. This is achieved by sampling the strains from Gaussian Processes (GPs). Each strain component is drawn from an independent multivariate normal distribution given by:
\begin{equation}
    \mathbf{X} \sim \mathcal{N}(\pmb{\mu}, \pmb{\Sigma})
\end{equation}
where $\mathbf{X}$ represents a vector containing the strain values at the different time steps, $\pmb{\mu}$ is the mean vector that specifies the expected value of strains, and $\pmb{\Sigma}$ is the covariance matrix. The covariance matrix $\pmb{\Sigma}$ describes the relationships between the samples in each of the components. The covariance function between two time steps $i$ and $j$ is given by:

\begin{equation}
    \Sigma_{ij} = k(x_i, x_j) = \sigma_f^2 \exp \left( - \frac{1}{2\ell^2} ||x_i - x_j||^2 \right)
\end{equation}
with $\sigma_f^2$ being the variance that determines the step size and $\ell$ being the length scale that controls the smoothness of the generated path. With increased variance $\sigma_f^2$ the strains are able to attain larger values, and with increased length scale the path becomes smoother. Values $\sigma_f^2 = 0.0001667$ and $\ell = 200$ are used in this work. A subset of the load paths generated by GPs is shown in Figure \ref{fig:gp}, with one path highlighted in red for clarity. We also show the strain paths for the highlighted loading path (see Figure \ref{fig:strains}) and the corresponding pair-wise stress-strain curves obtained from the full-order micromodel in Figure \ref{fig:strstr}.

% \begin{figure}[!h]
% \centering
% \includegraphics[scale=0.8]{Pictures/randomcurves.png}
% \caption{\centering Non-proportional loading path directions}
% \label{fig:gp}
% \end{figure}
\newpage
\section{Performance of PRNN with bulk model only}\label{prnn0}
\label{}
\begin{figure}[!h]
  \centering
  \includegraphics[scale=0.9]{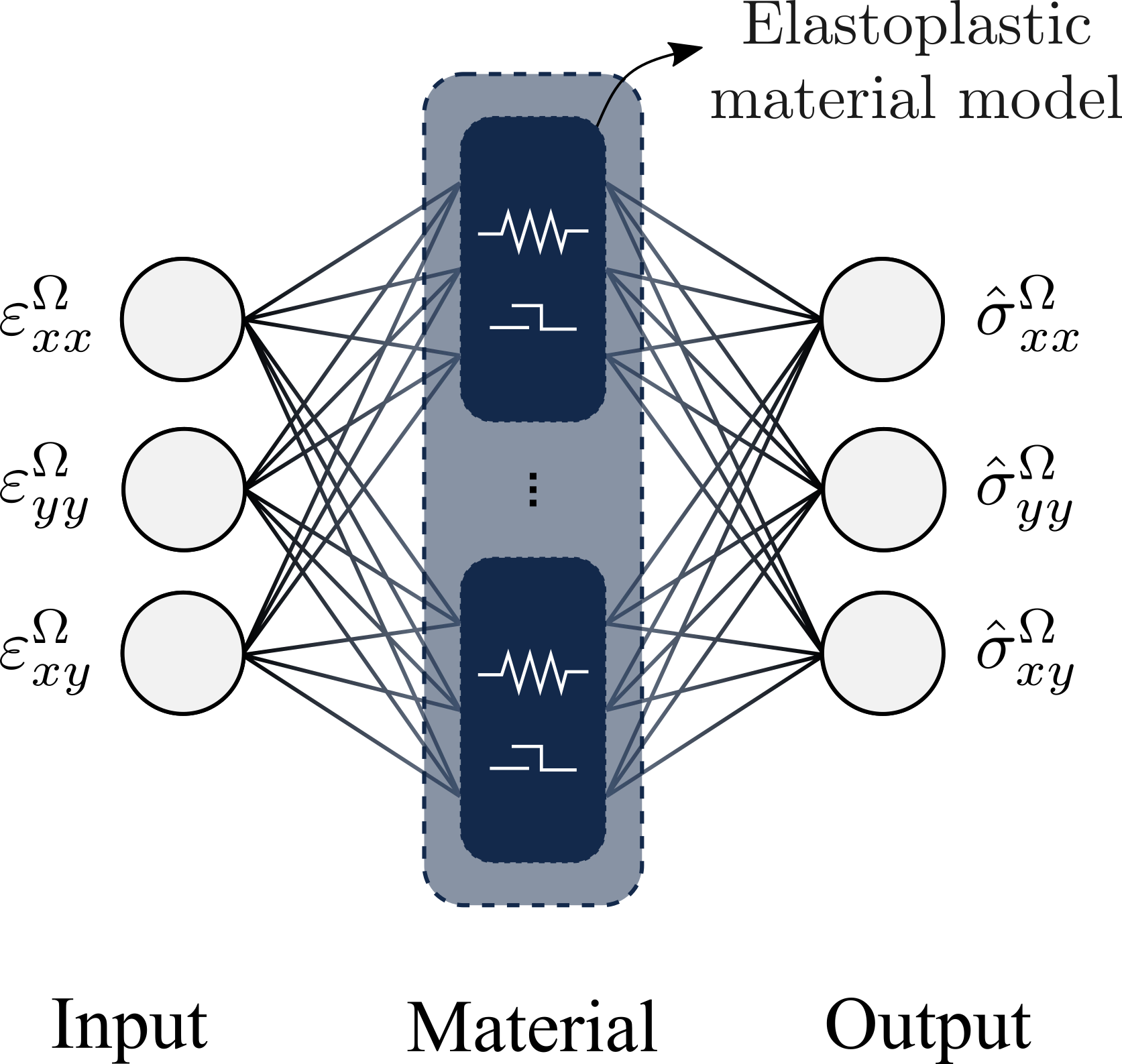}
  \caption{\centering PRNN with elastoplastic model only}\label{fig:prnn0}
\end{figure}

This section investigates whether the PRNN as proposed in \cite{MAIA2023115934} is able to capture stiffness degradation due to microscale damage. The architecture consists of one input layer, one material layer with bulk integration points only, and one output layer, as depicted in Figure \ref{fig:prnn0}. All bulk material points embed a J$_2$ plasticity model to convert 2D local strains to 2D local stresses.

\begin{figure}[!h]
\centering
\includegraphics[scale=0.6]{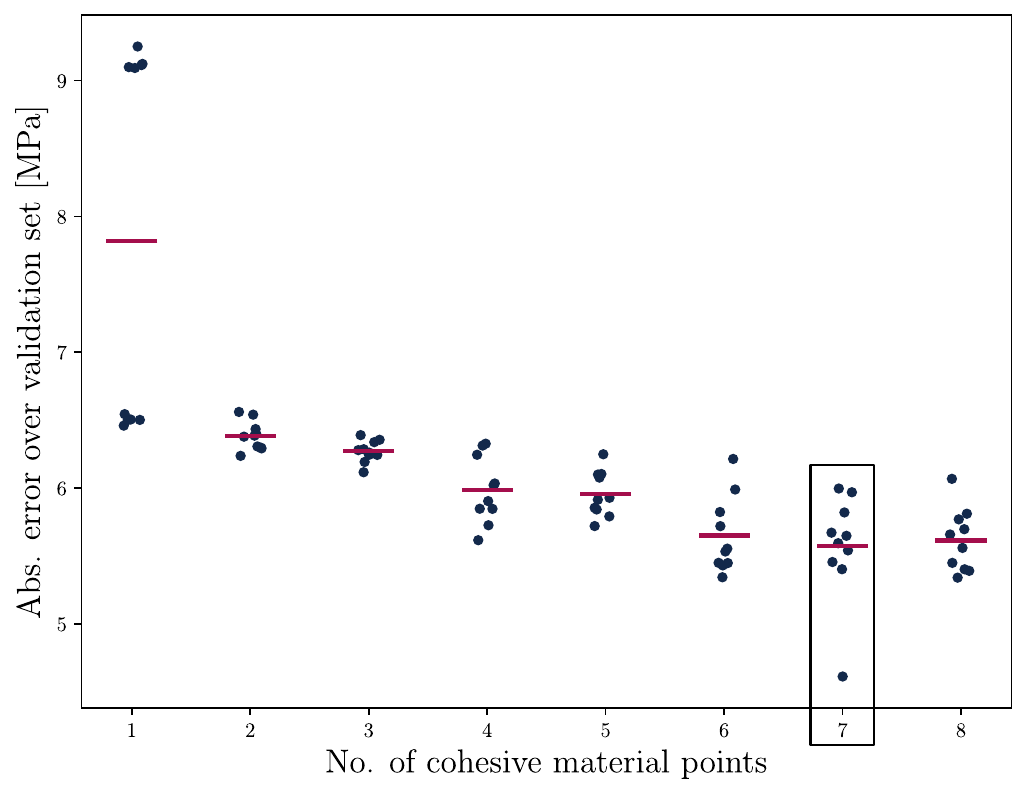}
\caption{\centering Validation error for PRNN trained on 18 fundamental curves with one cycle of unloading}
\label{fig:prnn0vali}
\end{figure}

In the original PRNN work \cite{MAIA2023115934}, the network could find a way to make elastoplastic material points reproduce linear elasticity by appropriately scaling encoder and decoder weights. In the following, we demonstrate how such an approach does not work for distributed damage. To highlight this inability to learn as clearly as possible, the networks here are trained and tested on the same curves. Specifically, the non-monotic, proportional dataset with one cycle of unloading in the 18 fundamental directions is used. Networks with different material sizes were trained by adding bulk points to the network until the mean value of the Mean Squared Errors (MSEs), no longer decreased with additional points. The best performing network with 7 fictitious material points and a training MSE of 4.61 MPa was selected for further examination. The training MSE across the different material layer sizes is shown in Figure \ref{fig:prnn0vali}, with the purple line representing the mean value for each material layer size.

% The validation set is identical to the training dataset, allowing for overfitting the training data.

\begin{figure}[h]
\centering
\begin{subfigure}[b]{0.49\textwidth}
\centering
\includegraphics[width=\linewidth]{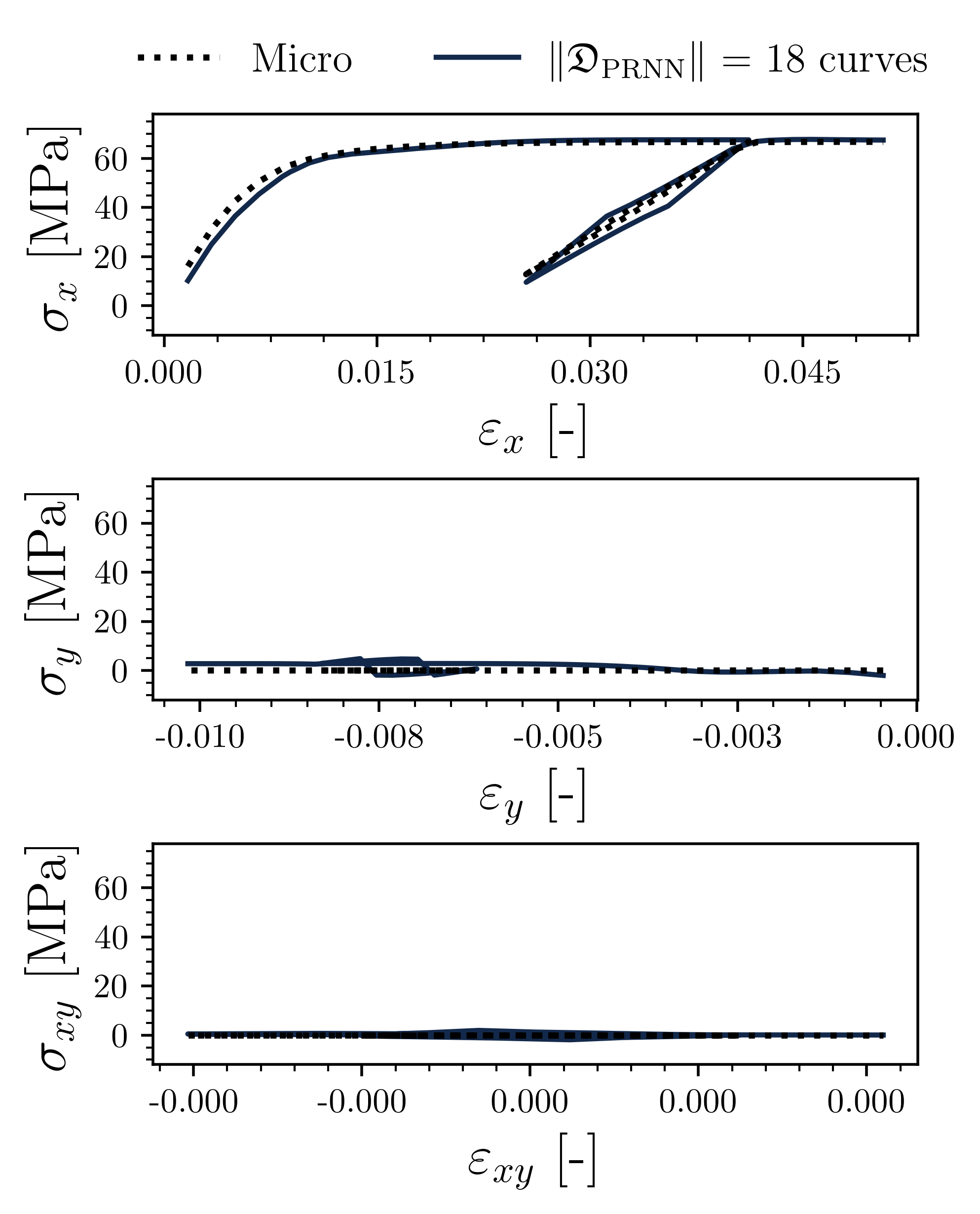}
\caption{\centering Uniaxial load case}
\label{fig:l06cr1}
\end{subfigure}
     \hfill
     \begin{subfigure}[b]{0.49\textwidth}
\centering
  \includegraphics[width=\linewidth]{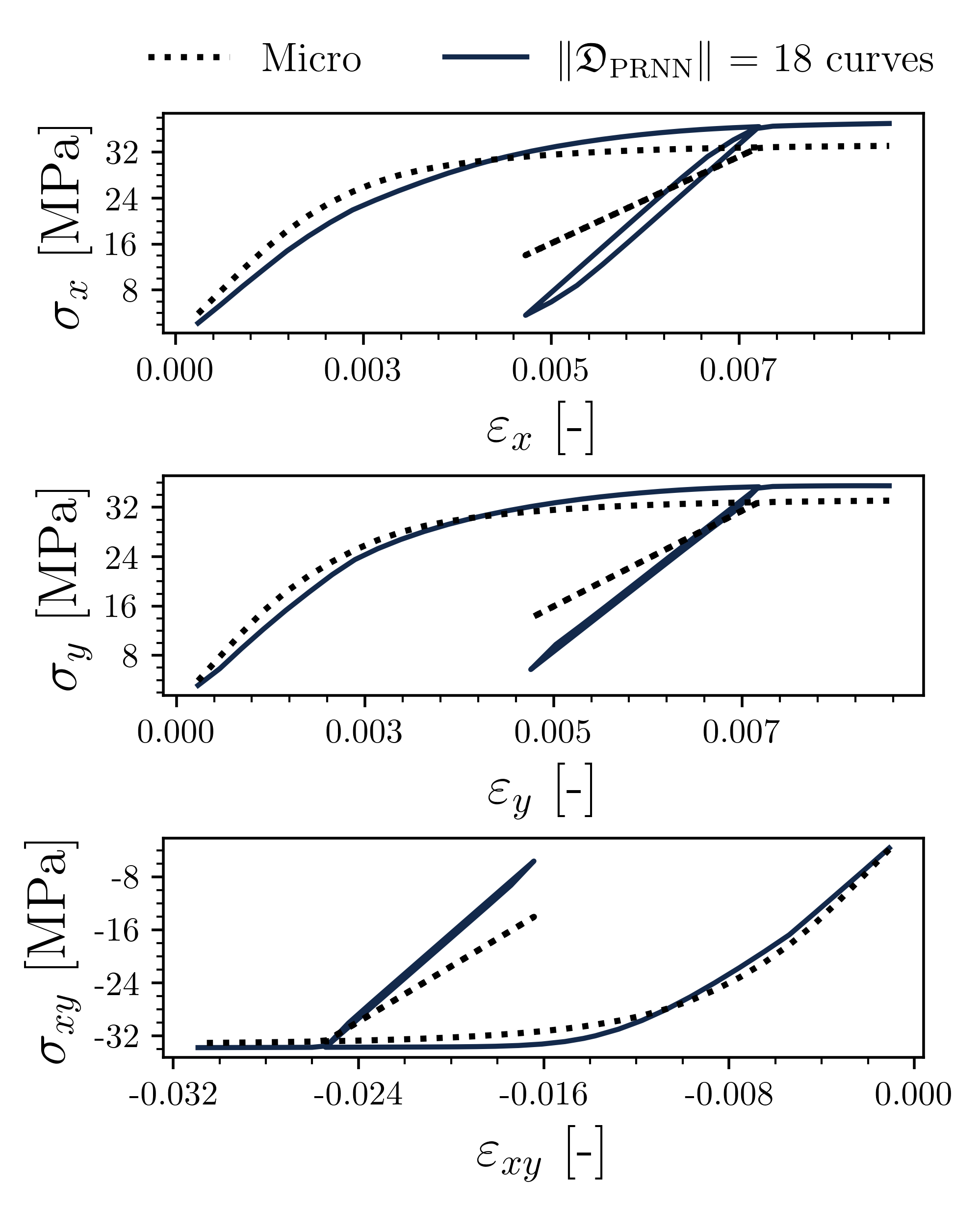}
  \caption{\centering Biaxial load case with shear}\label{fig:l06long}
\end{subfigure}
        \caption{\centering Prediction of PRNN trained with 18 fundamental curves with one cycle of unloading on training curves}
        \label{fig:monopred}
\end{figure}

The prediction of the network on two fundamental loading scenarios is shown in Figure \ref{fig:monopred}: uniaxial tension and biaxial tension with shear. The network provides a somewhat accurate prediction on the monotonic region of the curve, however, the model is unable to reproduce the unloading/reloading region. The PRNN starts to predict unloading with the initial, linear elastic stiffness following the assumptions embedded in the J$_2$ model, and predicts erratically afterwards. This highlights the limitation of the PRNN for describing stiffness degradation in its original design. The network encodes plasticity through the presence of a plasticity model in the material layer. This design gave the network a good bias in earlier work, when it could predict unloading behavior in plasticity without seeing it during training \cite{MAIA2023115934}. Here, however,the bias is too strong as it prevents the network from describing the stiffness loss that is present in the micromodel.

This observation is in line with the core idea of the PRNN to include a representation of all relevant physics by embedding the constitutive models from the micromodel in the network. This idea is violated by not including the cohesive zone model in the network. Therefore, the following sections of this work focus on the implementation of the cohesive zone model within the PRNN framework, along with evaluation of the proposed architectures.
\newpage
\section{Extending the network with cohesive material points}
\label{}
As shown in Section \ref{prnn0}, the physically recurrent neural network cannot accurately predict the effect of debonding at the fiber-matrix interface without including all sources of nonlinearity present in the RVE. Therefore, the cohesive zone model from the full-order micro-model has to be implemented in the PRNN as well. This section details the network configurations considered in this study for implementing the CZM within the PRNN framework.

\subsection{Cohesive points in the existing material layer}\label{cpem}

\begin{figure}[!h]
\minipage{0.49\textwidth}
  \centering
  \includegraphics[width=0.8\linewidth]{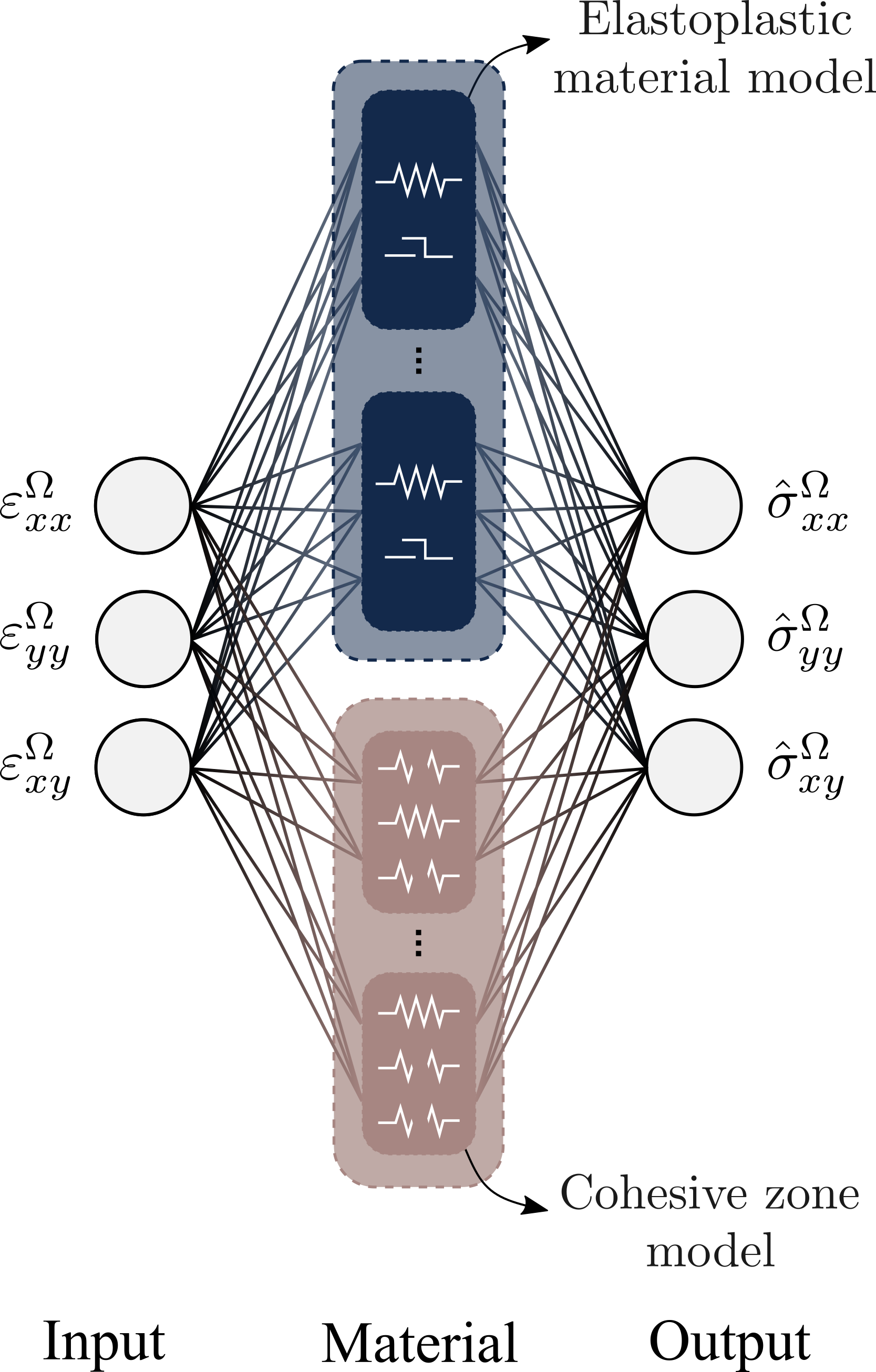}
  \caption{\centering PRNN$_1$ architecture: bulk and CZM in the same material layer}\label{fig:prnn1}
\endminipage\hfill
\minipage{0.49\textwidth}
  \centering
  \includegraphics[width=0.7\linewidth]{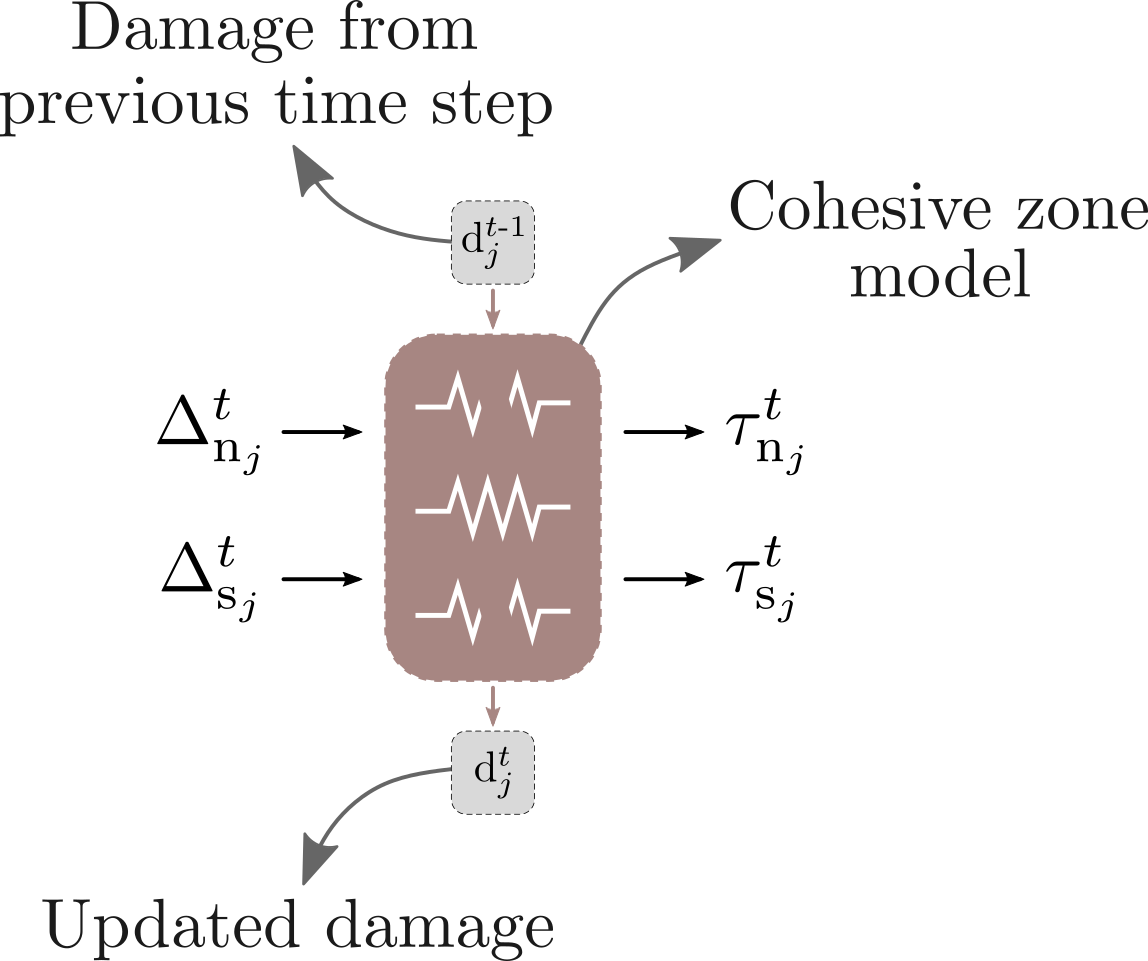}
  \caption{\centering Fictitious cohesive material point $j$}\label{fig:ficmatpt}
\endminipage\hfill
\end{figure}

The first design option retains the architecture proposed in \cite{MAIA2023115934} as much as possible. In this design, referred from now on as $\text{PRNN}_1$, there is one material layer containing bulk and cohesive fictitious points. The network is illustrated in Figure \ref{fig:prnn1} with bulk and cohesive points, represented in blue and pink, respectively. The cohesive points relate the local displacement jump vector, with normal and shear components, to a local traction vector, as illustrated in Figure \ref{fig:ficmatpt}. Similar to the bulk points, the cohesive points also store internal variables to account for history, in this case the damage variable defined as the ratio between dissipated energy and critical energy release rate \cite{TURON2018506}.

\subsection{Cohesive points in separate layer}

% This section investigates whether the PRNN with the CZM is able to capture micro-scale damage.

\begin{figure}[!h]
\centering
\includegraphics[scale=0.8]{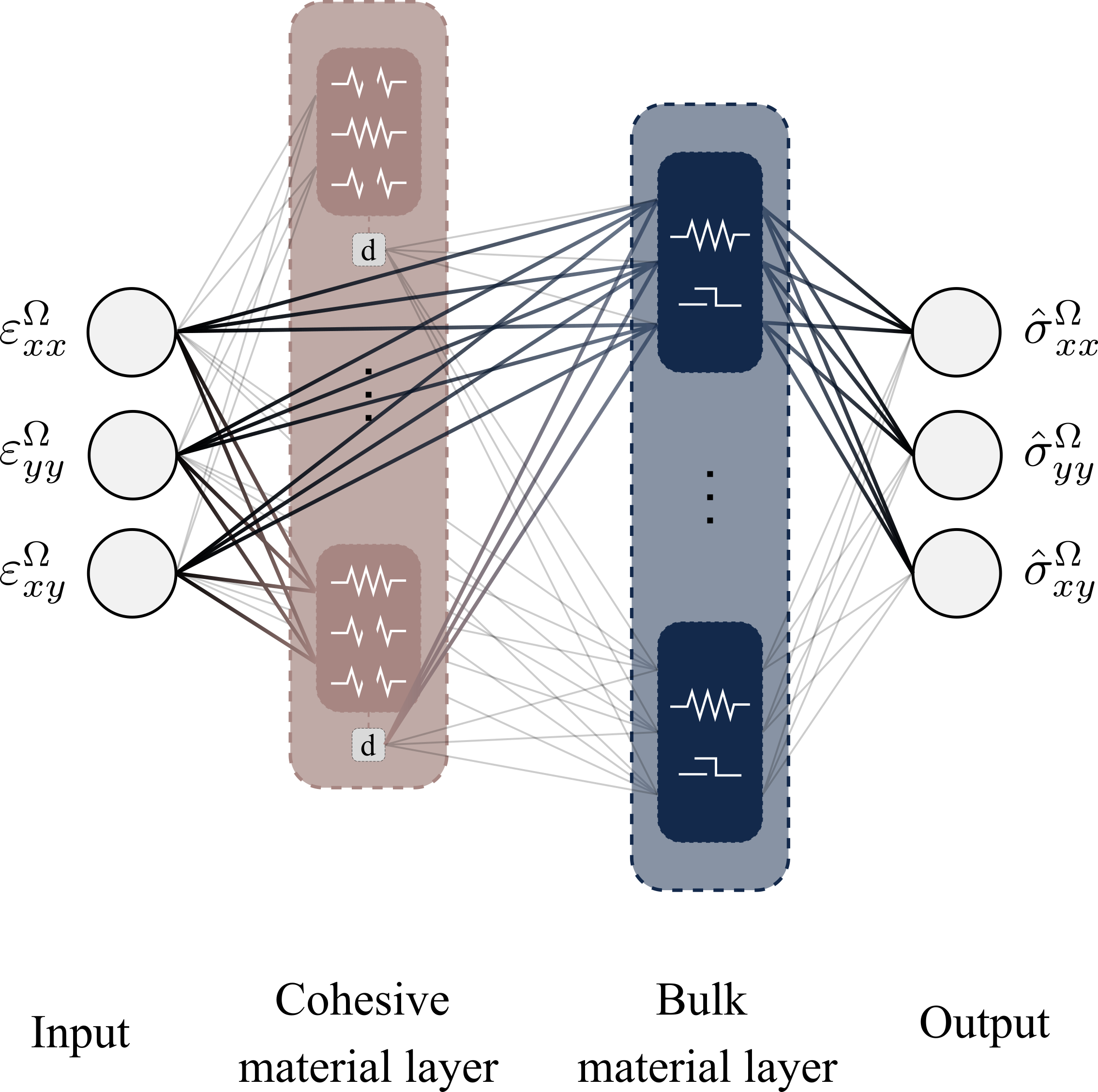}
\caption{\centering Novel architecture PRNN$_2$ and PRNN$_3$ with damage as input to bulk point}
\label{fig:dam}
\end{figure}

Instead of having both types of models in the same layer, we also investigate architectures with two material layers: a cohesive and a bulk material layer, each with embedded models, as illustrated in Figure \ref{fig:dam}. The two architectures considered here consist of one input layer that receives macroscopic strain, two material layers containing the nonlinear models, and one output layer yielding the macro-scale stress predictions. The state variables of the cohesive points are densely connected to the bulk points together with the macroscopic strains. To illustrate the connectivity of the layers, Figure \ref{fig:dam} highlights how one cohesive and one bulk are linked to each other, as well as to the input and output.

Rather than using the output traction values of the cohesive points, damage is used as input to the bulk points. Damage, the internal variable stored in the cohesive points, either increases or remains the same in case of unloading, providing a more monotonic influence on the overall response. This damage variable modifies the local strain value received by the bulk points, resulting in adjusted local stress values for the same level of macroscopic strain. The irreversibility of damage gives rise to a decrease in stiffness during unloading. This design ensures that only the bulk points contribute directly to the stress homogenization procedure, unlike in the architecture described in Section \ref{cpem}, where tractions rising from the cohesive points are directly connected to the output through the decoder layer. This is more consistent with the homogenization procedure in FE$^2$, where cohesive tractions do not contribute directly to the macroscopic stress (cf. Eq \eqref{eq:homogenization}).

Two ways of connecting the damage variable from the cohesive points to the local strain at the bulk points are considered. The first method follows a more conventional approach, which involves densely connecting the damage values to the bulk points:

\begin{equation}
\pmb{\varepsilon} = \mathbf{w}_\text{d} \cdot \mathbf{d} + \mathbf{w}_{\varepsilon \text{b}} \cdot \pmb{\varepsilon}^{\Omega}
\end{equation}
where $\mathbf{w}_\text{d}$ and $\mathbf{w}_{\varepsilon \text{b}}$ are the weight matrices connecting the damage values $\mathbf{d}$ from all the cohesive points and the macroscale strain, respectively, to the local strain values of the bulk points. The network with this approach will be referred to as $\text{PRNN}_2$ from now on.

In the second method, referred to as $\text{PRNN}_3$ from now on, the damage variables are used to modify the amplitude of the local strain input to the bulk points. This is achieved by multiplying the local strain input piece-wise by the term $\phi_{\text{sp}}(\mathbf{1}+\mathbf{w}_\text{d} \cdot \mathbf{d}))$, which is forced to obtain positive values by applying a SoftPlus activation function ($\phi_{\text{sp}}(\cdot)$):

\begin{equation}
\pmb{\varepsilon} = \phi_{\text{sp}}(\mathbf{1}+\mathbf{w}_\text{d} \cdot \mathbf{d}) *(\mathbf{w}_{\varepsilon \text{b}} \cdot \pmb{\varepsilon}^{\Omega})
\end{equation}

\section{Performance of PRNN with cohesive model}
\label{}
In this section, we assess the performance of the PRNNs with the architectures proposed in Section 5. The model selection process is presented by analyzing their performance across different training set and material layer sizes. We evaluate and compare the PRNNs' ability to accurately capture the homogenized response, taking into account microscale damage, under various loading scenarios.

\subsection{Model selection}

\begin{figure}[!h]
\centering
\includegraphics[width=0.48\linewidth]{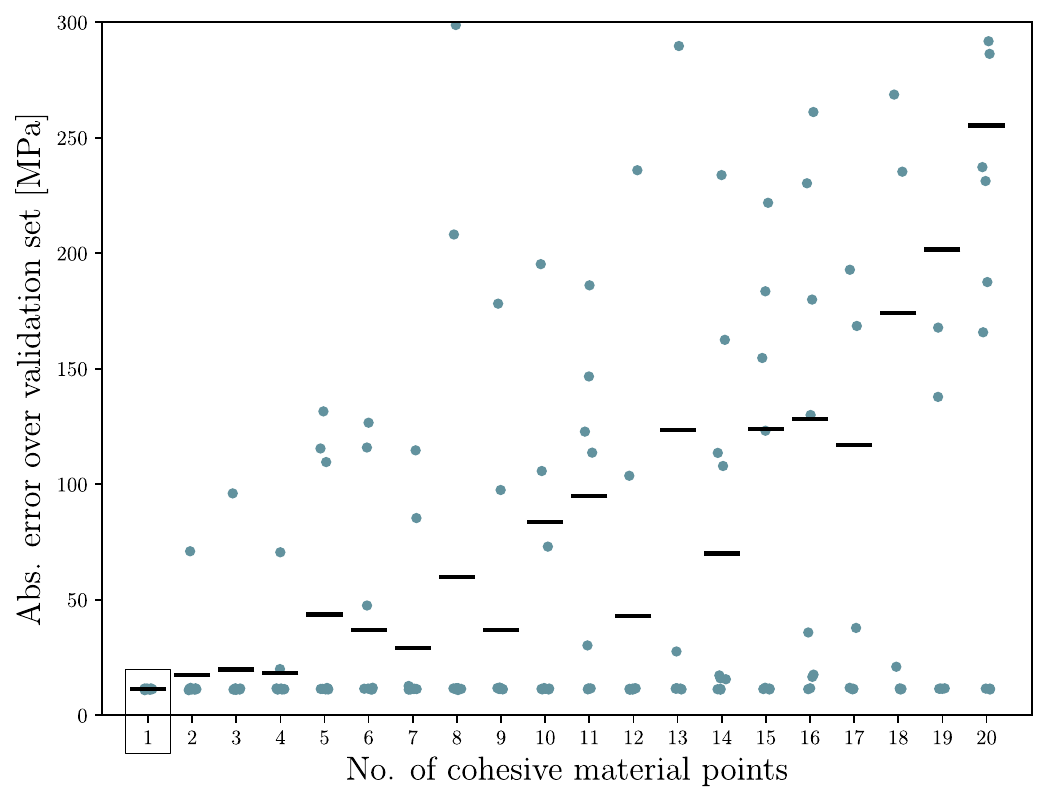}
\caption{\centering Validation error for $\text{PRNN}_1$ trained on 192 GP curves}
\label{fig:gp72}
\end{figure}

\begin{figure}[!h]
\minipage{0.48\textwidth}
  \centering
  \includegraphics[width=1.\linewidth]{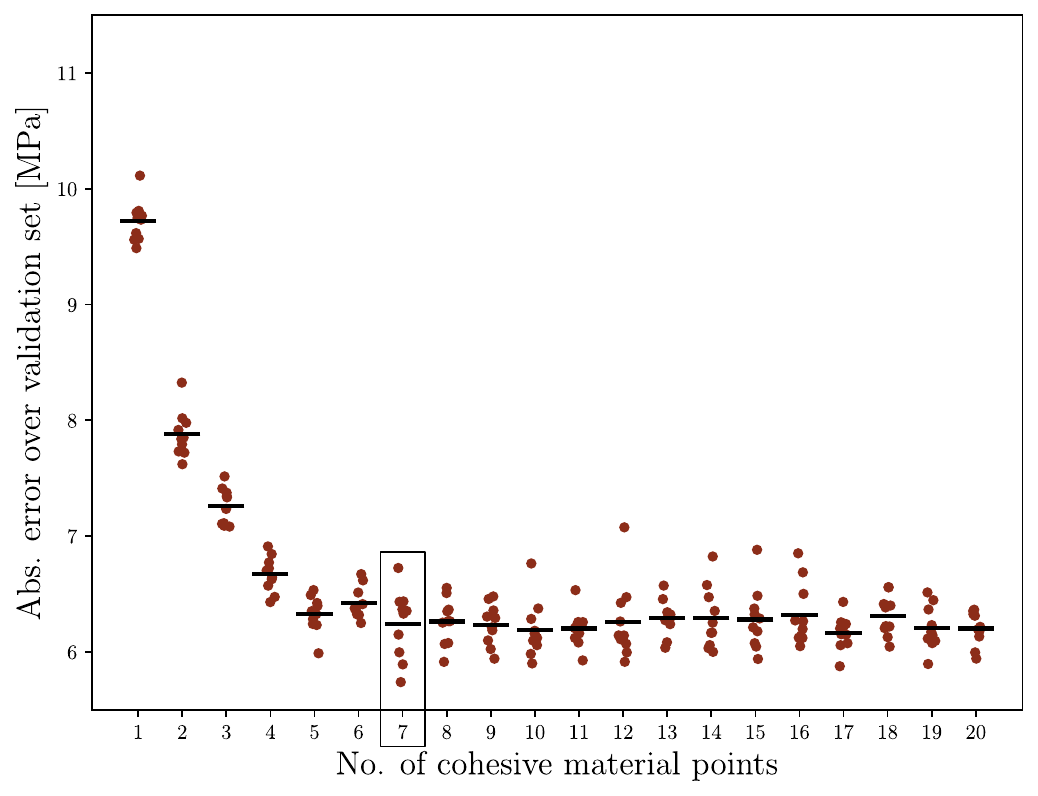}
\caption{\centering Validation error for $\text{PRNN}_2$ trained on 192 GP curves}
\label{fig:6v}
\endminipage\hfill
\minipage{0.48\textwidth}
  \centering
  \includegraphics[width=1.\linewidth]{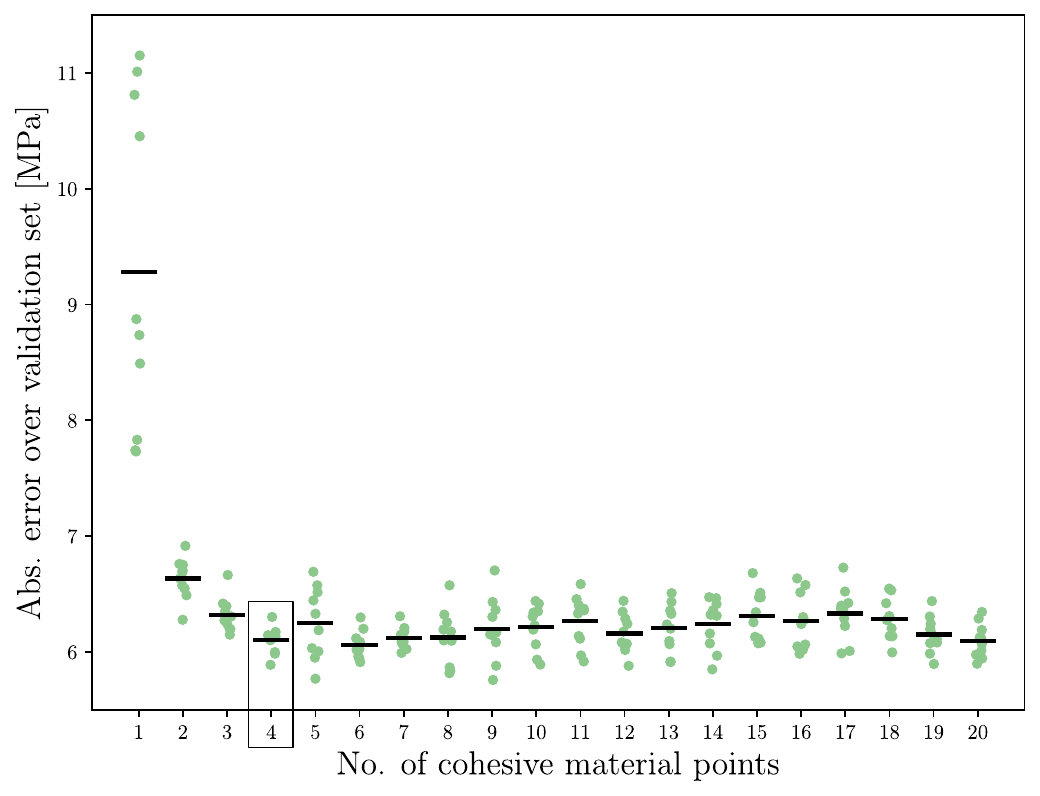}
\caption{\centering Validation error for $\text{PRNN}_3$ trained on 192 GP curves}
\label{fig:8v}
\endminipage\hfill
\end{figure}

First we perform model selection for the size of the material layer. For that purpose, networks were trained on 192 GP-based curves (non-monotonic and non-proportional loading) with varying numbers of bulk and cohesive points. The ratio of bulk to cohesive points is kept constant and equal to 4, mirroring the ratio of matrix to cohesive elements in the RVE. The size of the material layer ranges from a minimum configuration of four bulk and one cohesive points to 80 bulk and 20 cohesive points. Figures \ref{fig:gp72}-\ref{fig:8v} show the MSE on a validation set with 200 GP-based curves for 10 different initializations in each material layer size, for all three architectures considered in this work. For $\text{PRNN}_1$, networks with 4 bulk points and 1 cohesive point are selected, while for $\text{PRNN}_2$ a combination with 28 bulk points and 7 cohesive points is needed. Finally, for the $\text{PRNN}_3$, 44 bulk and 11 cohesive points are selected. The observation that the validation error increases for increasing network sizes of PRNN$_1$ points at the tendency of this architecture to overfit.  

\begin{figure}[!h]
\centering
\includegraphics[scale=0.45]{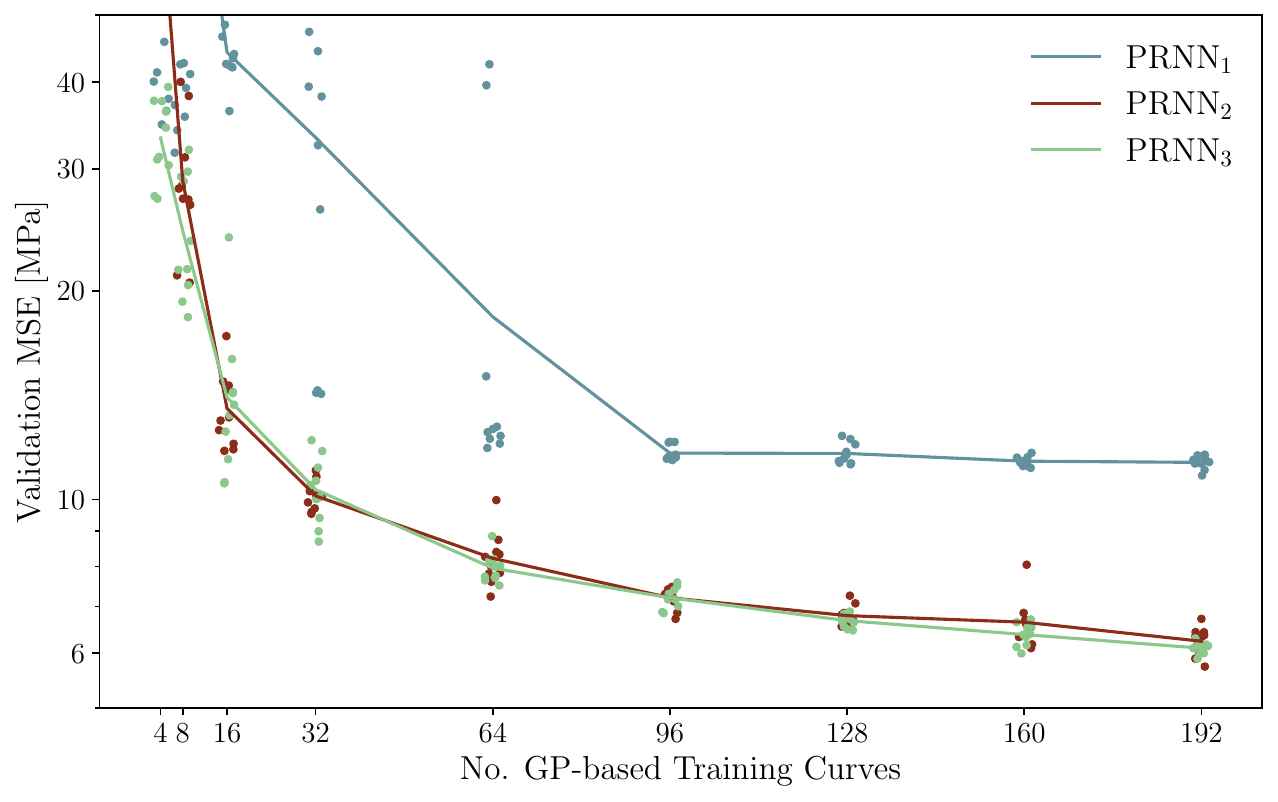}
\caption{\centering Validation MSE for the PRNNs considered across various training data sizes}
\label{fig:gpvalitr}
\end{figure}

\newpage
The selected networks were trained on different training set sizes, ranging from 4 to 192 GP-based curves (non-monotonic and non-proportional loading). Figure \ref{fig:gpvalitr} displays the MSE on a validation set with 200 GP-based curves across the various training data sizes for the selected material layer sizes of the three architectures considered in this work. The solid lines in the figure represent the mean MSE values for each PRNN at different training data sizes. Based on Figure \ref{fig:gpvalitr}, a training set size of 96 paths was selected for the first architecture ($\text{PRNN}_1$) and 192 for the latest two architectures ($\text{PRNN}_2$ and $\text{PRNN}_3$). The lowest validation error corresponding to these networks is 11.40MPa, 5.74MPa, and 5.89MPa, respectively.

%, while fig \textcolor{red}{\ref{fig:gp68}} shows the comparison for 6 and 8 for a larger number of training curves. 
%It is important to note that the validation error decreases significantly when the cohesive points are implemented in a seperate layer from the bulk points, quantitatively indicating that the modified architecture has an improved prediction accuracy. 

% \begin{figure}[!h]
% \centering
% \includegraphics[scale=0.6]{elsarticle/Pictures/boxcomparezANN6v_4ii_192gp.pdf}
% \caption{\centering Validation error for $\text{PRNN}_2$ trained on 192 GP curves}
% \label{fig:gp72}
% \end{figure}

% \begin{figure}[!h]
% \centering
% \includegraphics[scale=0.6]{elsarticle/Pictures/boxcomparezANN8v_4ii_192gp.pdf}
% \caption{\centering Validation error for $\text{PRNN}_3$ trained on 192 GP curves}
% \label{fig:gp72}
% \end{figure}

% \begin{figure}[!h]
% \centering
% \includegraphics[scale=0.6]{elsarticle/Pictures/ANN68.pdf}
% \caption{\centering Lowest validation MSE for different material layer sizes across various training data sizes for the modified architecture}
% \label{fig:gp68}
% \end{figure}

% A training data size of 192 was selected, as the average value of the lowest MSE for different material layer sizes does not decrease significantly when further increasing the training size. Then, based on the MSE obtained on the validation set, a network with a combination of 4 bulk and 1 cohesive points was selected with a lowest MSE of \textcolor{red}{?} MPa. 

% \newpage
\subsection{Predicting micro-scale damage}

The selected networks were tested on two different datasets: one containing 54 curves from the random, non-monotonic, and proportional dataset with one cycle of unloading, and another with 54 curves from the non-monotonic, non-proportional dataset (different from the ones used for training and validation but of the same type). The average MSE on these two test sets are presented in Table \ref{tab} for each network.

\begin{table}[h]
\centering
\begin{tabular}{lllll}
\cline{1-4}
\multicolumn{1}{|l|}{Test set} & \multicolumn{1}{l|}{$\text{PRNN}_1$} & \multicolumn{1}{l|}{$\text{PRNN}_2$} & \multicolumn{1}{l|}{$\text{PRNN}_3$} &  \\ \cline{1-4}
\multicolumn{1}{|l|}{Non-monotonic, proportional} & \multicolumn{1}{c|}{7.45} & \multicolumn{1}{c|}{5.56} & \multicolumn{1}{c|}{3.40} &  \\ \cline{1-4}
\multicolumn{1}{|l|}{Non-monotonic, non-proportional} & \multicolumn{1}{c|}{11.63} & \multicolumn{1}{c|}{5.72} & \multicolumn{1}{c|}{6.03} & \\ \cline{1-4}
\end{tabular}
\caption{Average MSE values for the two test datasets, in MPa}
\label{tab}
\end{table}

\begin{figure}[!h]
  \centering
  \begin{subfigure}[b]{0.45\textwidth}
    \centering
    \includegraphics[width=\textwidth]{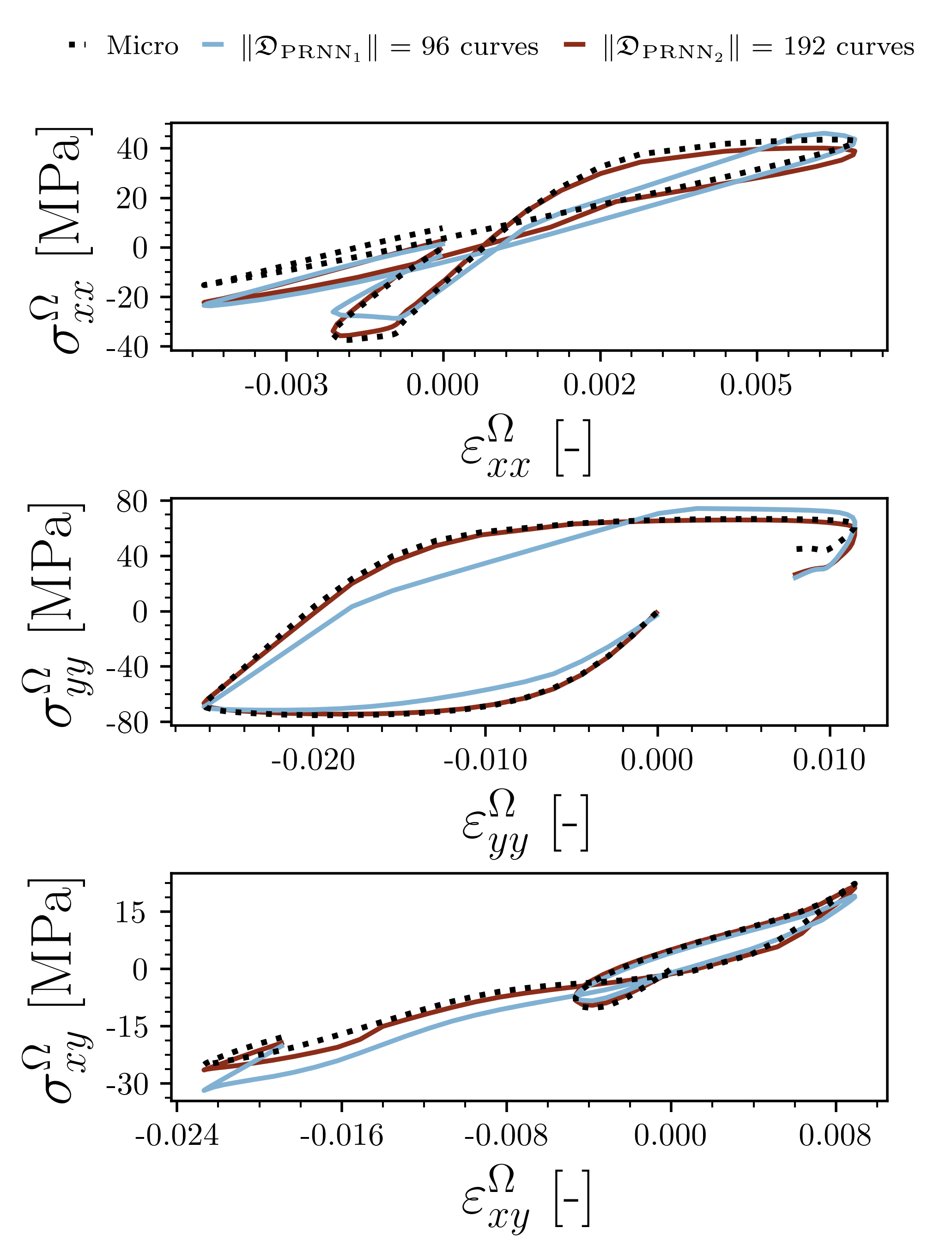}
    \caption{\centering Prediction of $\text{PRNN}_1$ and $\text{PRNN}_2$ on a representative GP curve}
    \label{fig:v6vgp}
  \end{subfigure}
  \hfill
  \begin{subfigure}[b]{0.45\textwidth}
    \centering
    \includegraphics[width=\textwidth]{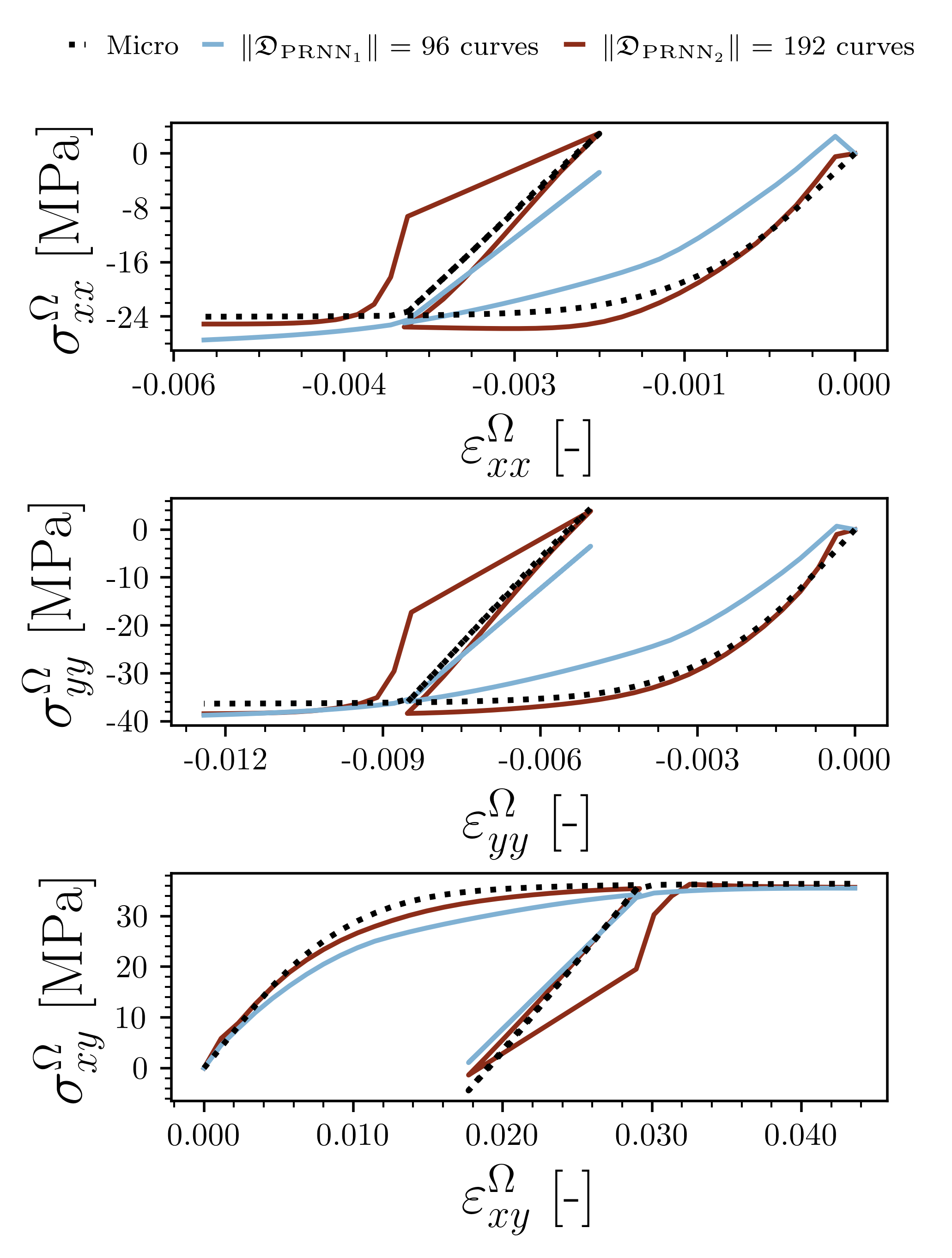}
    \caption{\centering Prediction of $\text{PRNN}_1$ and $\text{PRNN}_2$ on a representative curve with one cycle of unloading}
    \label{fig:v6vnm}
  \end{subfigure}
\caption{\centering Prediction of $\text{PRNN}_1$ and $\text{PRNN}_2$ on representative test curves}
  \label{fig:v6v}
\end{figure}

First, we illustrate the performance of $\text{PRNN}_1$ and $\text{PRNN}_2$ using a representative curve from each test set in Figure \ref{fig:v6v}, with MSE close to the average MSE shown in Table \ref{tab}. Note how the predictions of $\text{PRNN}_1$ on GP-based curves follow the overall trend but with significantly less accuracy compared to $\text{PRNN}_2$ (Figure \ref{fig:v6vgp}). This observation aligns well with the results shown in Figure \ref{fig:gpvalitr}, where the validation set is also comprised of GP-based curves, emphasizing the significant decrease in validation error when the cohesive points are implemented in a separate layer from the bulk points. The difference between $\text{PRNN}_1$ and $\text{PRNN}_2$ predictions becomes more pronounced when tested on the non-monotonic, proportional dataset. As shown in Figure \ref{fig:v6vnm}, $\text{PRNN}_1$ predicts poorly. The network not only fails capturing the decrease in stiffness during unloading but also loses accuracy in the monotonic part. It is observed that with $\text{PRNN}_1$ there is a preference towards networks with fewer cohesive points. Moreover, small weights connect the normal component of these cohesive points with the output, which is likely due to the large traction values output from the cohesive points in compression. These factors indicate that the network avoids utilizing the cohesive points implemented in the material layer, resulting in unloading with the initial linear stiffness.

% The issue also points to an important inconsistency between the architecture of the network and the underlying physics in the full-order solution, mentioned in section \ref{cpem}. In the $\text{FE}^2$ method, introduced in Section \ref{fenn}, the micro-scale tractions computed by interface elements interact with the bulk points in the equilibrium equations of the RVE. These micro-scale tractions, associated to interface elements, do not contribute to the homogenized macro-scale stresses in the RVE; only the micro-scale stresses, associated with elements with volume, from the bulk points do. However, in $\text{PRNN}_1$, the micro-scale tractions are connected to the output with a dense layer. This arrangement does not adequately represent the interaction of the cohesive and bulk points in the equilibrium problem, nor does it capture the correct (stress) contributions to the homogenization procedure.

\begin{figure}[!h]
  \centering
  \begin{subfigure}[b]{0.6\textwidth}
    \centering
    \includegraphics[width=\textwidth]{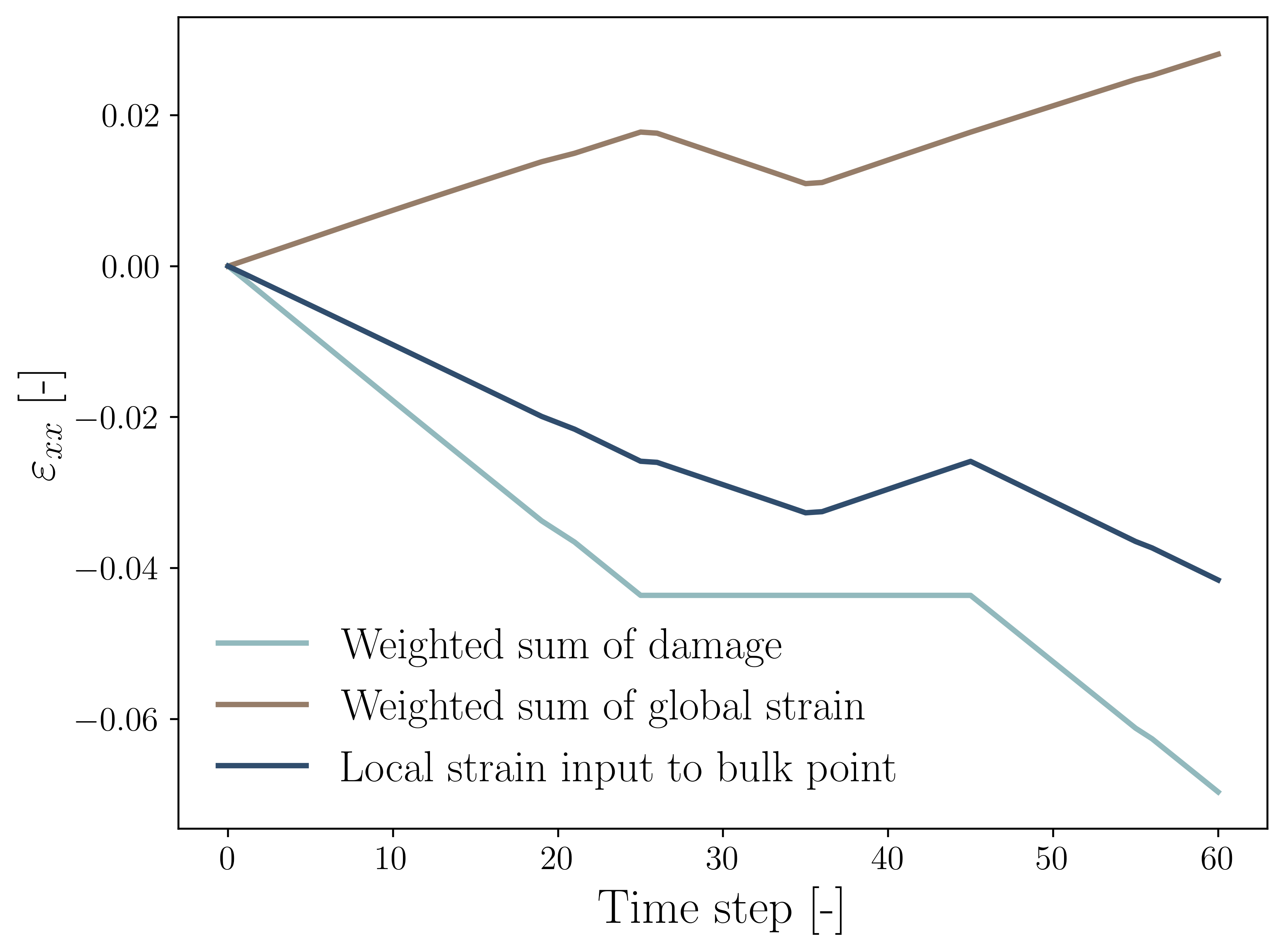}
    \caption{\centering Input to the y-component of the bulk point}
    \label{fig:m19}
  \end{subfigure}
  \vfill
  \begin{subfigure}[b]{1.\textwidth}
    \centering
    \includegraphics[width=\textwidth]{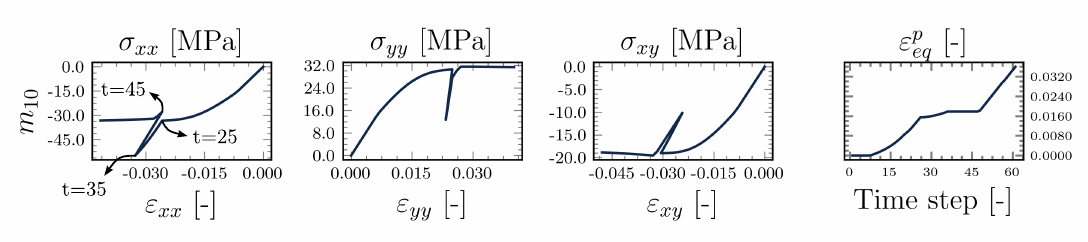}
    \caption{\centering Response of the bulk point}
    \label{fig:m19bulk}
  \end{subfigure}
\caption{\centering Behavior of one of the bulk points of $\text{PRNN}_2$ when predicting on a representative curve from the non-monotonic, propotional test set}
  \label{fig:m199}
\end{figure}

Besides the improved accuracy on the test sets, as shown in Table \ref{tab}, $\text{PRNN}_2$ shows another advantage. Unloading occurs with a different slope than the initial linear phase (Figure \ref{fig:v6vnm}), indicating that the network is able to account for the effect of microscale damage. However, a new problem arises: reloading follows a different path than unloading. This phenomenon is explained by closely examining the input to one of the bulk points. Figure \ref{fig:m19} illustrates the $\varepsilon_{yy}$ component of a particular fictitious bulk material point and its two contributions, one that follows directly from the macroscopic strain ($\mathbf{w}_{\varepsilon \text{b}} \cdot \pmb{\varepsilon}^{\Omega}$) and the other that follows from the damage variables from all the cohesive points ($\mathbf{w}_\text{d} \cdot \mathbf{d}$). Globally, the micromodel is unloading, a trend represented by the weighted sum of the global strain values. Meanwhile, the weighted sum of damage is larger than the weighted sum of the global strain and is with opposite sign. Therefore, the final sum used as input to the bulk point not only changes signs but also prevents the point to follow the global unloading trend (from t = 25 to t = 35). Instead, the bulk point continues to load while the macroscopic strain is subjected to unloading and only starts unloading once the macroscopic strain is at the reloading branch (from t = 35 to t = 45). This mismatch between the loading phases leads to further evolution of plastic strain during the macroscopic unloading, as shown in Figure (\ref{fig:m19bulk}), and causes the undesired change in slopes during unloading and reloading.

\begin{figure}[!h]
  \centering
  \begin{subfigure}[b]{0.44\textwidth}
    \centering
    \includegraphics[width=\textwidth]{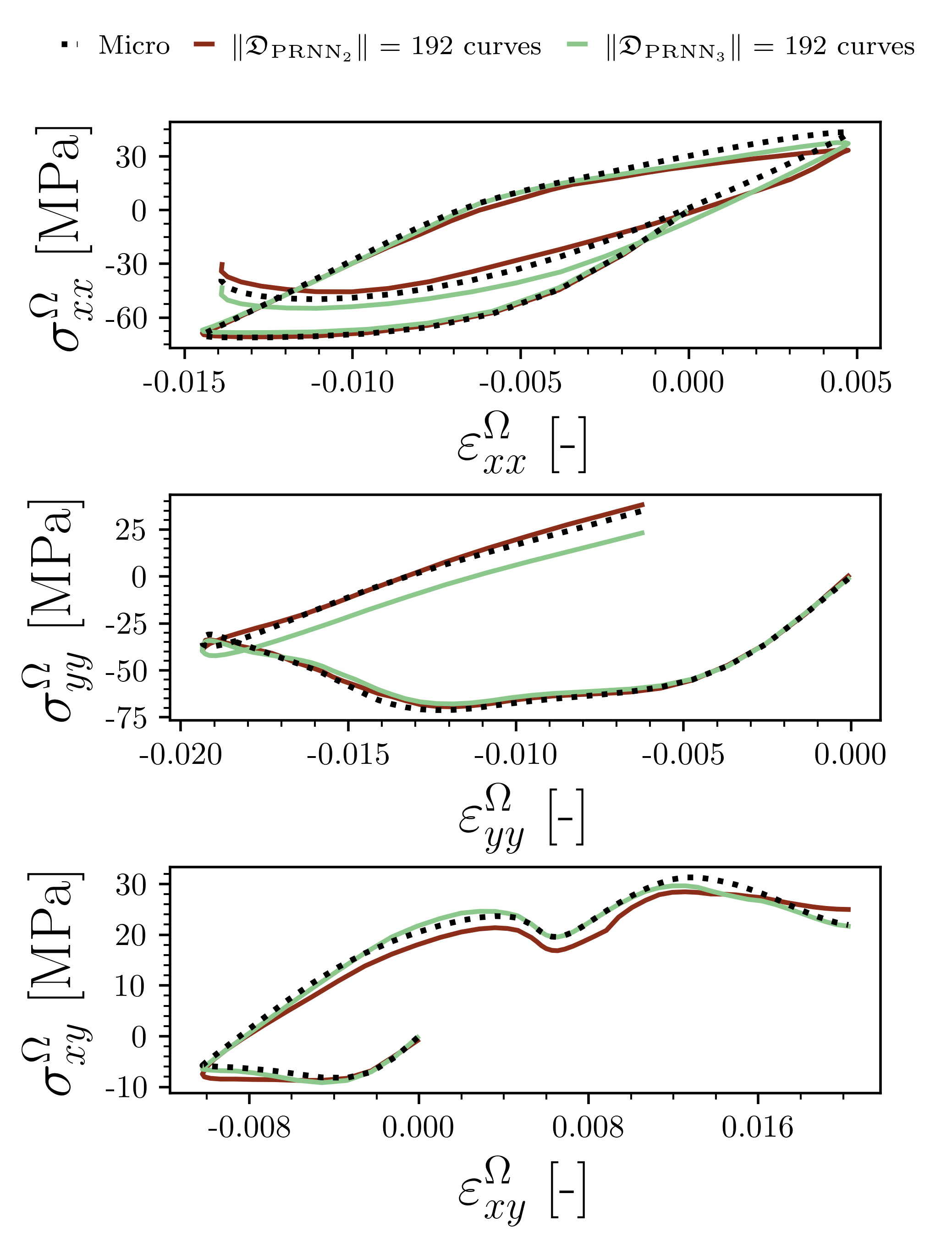}
    \caption{\centering Prediction of $\text{PRNN}_2$ and $\text{PRNN}_3$ on a representative GP curve}
    \label{fig:6v8vgp}
  \end{subfigure}
  \hfill
  \begin{subfigure}[b]{0.44\textwidth}
    \centering
    \includegraphics[width=\textwidth]{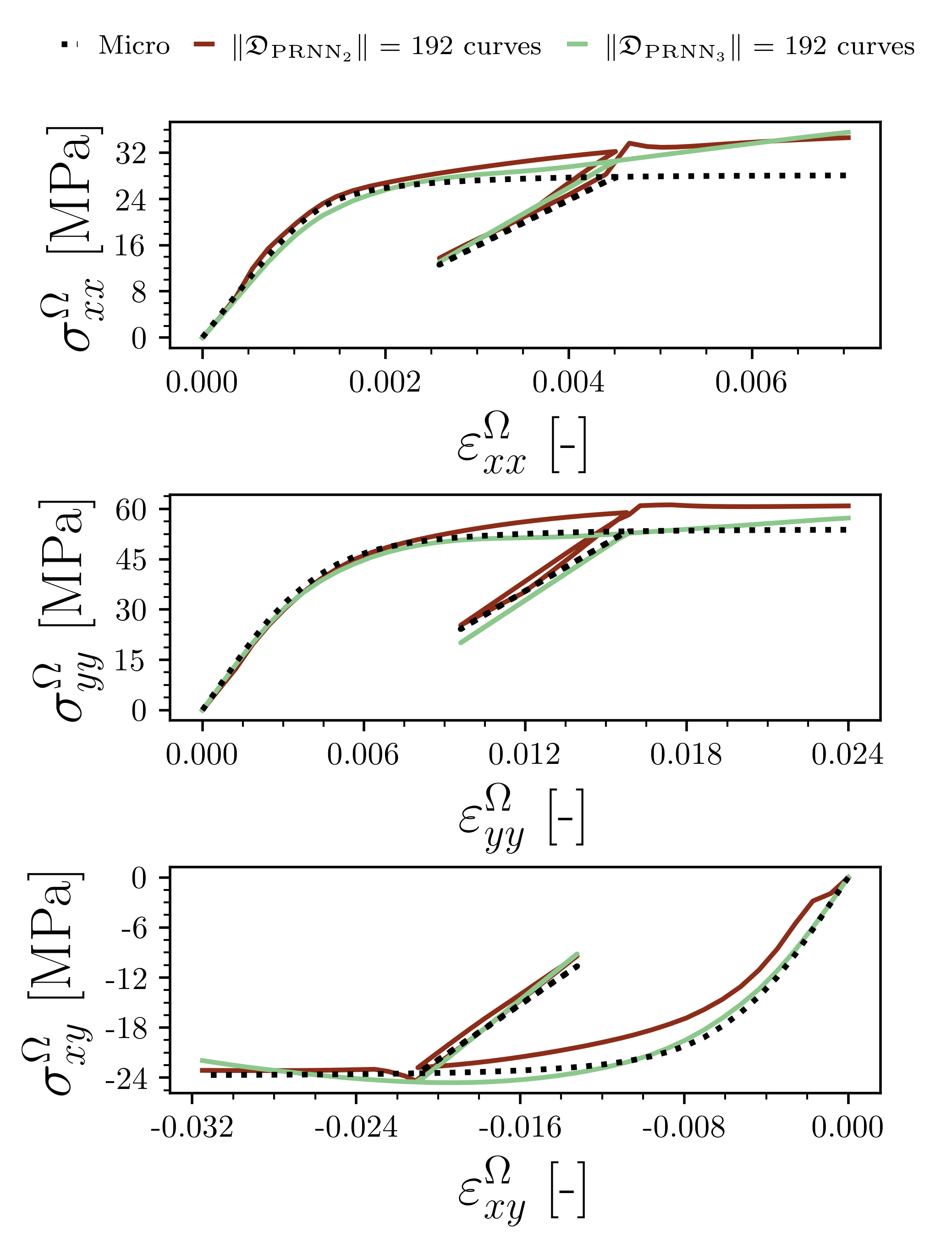}
    \caption{\centering Prediction of $\text{PRNN}_2$ and $\text{PRNN}_3$ on a representative curve with one cycle of unloading}
    \label{fig:6v8vnm}
  \end{subfigure}
\caption{\centering Prediction of $\text{PRNN}_2$ and $\text{PRNN}_3$ on representative test curves}
  \label{fig:6v8v}
\end{figure}

When damage is used as an amplifier rather than being simply densely connected to the bulk points, no significant difference in validation errors was found (Figure \ref{fig:gpvalitr}). This is reflected on the prediction of $\text{PRNN}_2$ and $\text{PRNN}_3$ on a representative GP-based curve in Figure \ref{fig:6v8vgp}. However, robustness improves significantly when predicting on curves from the non-monotonic, proportional set. Figure \ref{fig:6v8vnm} clearly shows that unloading/reloading now takes place along the same path and with a different slope than the initial linear phase, effectively capturing the effect of microscale damage.

\begin{figure}[!h]
\centering
\includegraphics[scale=0.73]{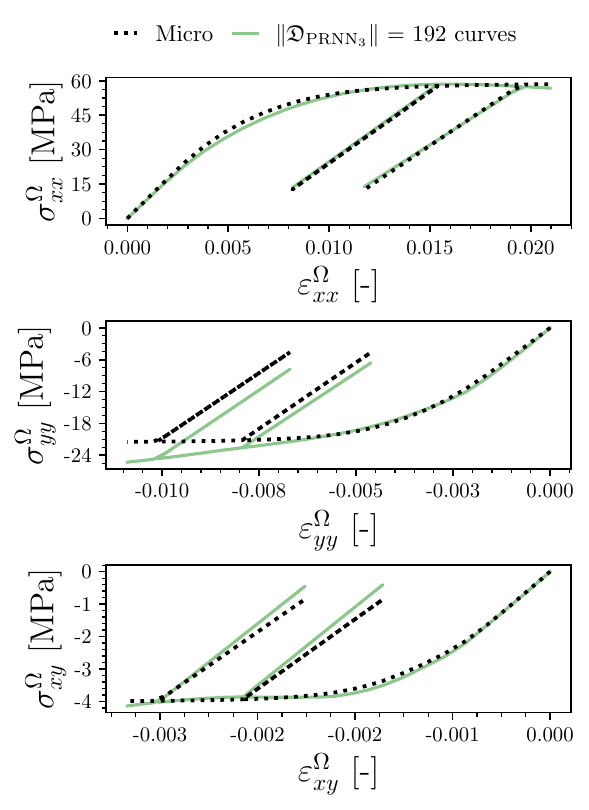}
\caption{\centering Prediction of $\text{PRNN}_3$ on a representative curve with two cycles of unloading}
\label{fig:twooo}
\end{figure}

To further demonstrate the network's predictive capabilities, Figure \ref{fig:twooo} shows the prediction of $\text{PRNN}_3$ on a curve from a test set containing non-monotonic, proportional curves with two cycles of unloading. The evolution of damage over time is evident as the slope of the unloading-reloading phase gradually decreases as the loading continues.

% \begin{figure}[!h]
%   \centering
%   \begin{subfigure}[b]{0.47\textwidth}
%     \centering
%     \includegraphics[width=\textwidth]{elsarticle/Pictures/ANN8vmod_4812_192gp_52_1.png}
%     \caption{\centering Prediction of network on a representative curve with one unloading}
%     \label{fig:23gp72good}
%   \end{subfigure}
%   \hfill
%   \begin{subfigure}[b]{0.47\textwidth}
%     \centering
%     \includegraphics[width=\textwidth]{elsarticle/Pictures/ANN8vmod_4812_192gp_nm2_9_1.png}
%     \caption{\centering Prediction of network on a representative curve with two unloadings}
%     \label{fig:23gp72avg}
%   \end{subfigure}

%   \caption{\centering Prediction of network with 2 bulk and 3 cohesive points trained on 72 GP curves on test curves with one cycle of unloading}
%   \label{fig:23gp72}
% \end{figure}
\newpage
\section{Conclusions}
\label{}
In this paper, we propose an extension to a recently proposed surrogate model, namely the Physically Recurrent Neural Network (PRNN), to account for the complex combination of plasticity and microscale damage. The PRNN's excellent ability to predict elastoplastic behavior motivated this study into its use as a surrogate model in a more challenging context where both plasticity and damage are present. Constitutive relations from the full-order micromodel are directly implemented into the hidden layer of the PRNN, creating a direct link to the micromodel. Path-dependency naturally arises from the material models in the network, resulting in accurate predictions with a significantly smaller training dataset compared to networks without physical interpretation.

As a first step to use the PRNN framework for microscale damage, a preliminary study was conducted using the network in its original form with bulk material points only. It was demonstrated that the original PRNN could not describe stiffness degradation, even when trying to overfit on a small set of training curves. These results align well with the general guideline in \cite{MAIA2023115934} that all types of nonlinearities present in the RVE need to be included in the network, justifying the need for an extended PRNN architecture that integrates a cohesive zone model.

Next, three architectures of the PRNN with bulk and cohesive points were proposed. In the first design, $\text{PRNN}_1$, cohesive points with the CZM were incorporated into the same material layer of the PRNN as the plasticity model. Two material layers were used in the second and third design with the cohesive points implemented in a separate cohesive layer from the bulk material layer. Together with the global strain, the internal variable of the cohesive points, damage, was then used as input to the bulk points. This connection was defined in two different ways, either in a conventional way with a dense connection ($\text{PRNN}_2$) or by using the damage as an amplifier to the local strain of the bulk points ($\text{PRNN}_3$). %Additionally, to account for the highly nonlinear relation between macroscale strains and local displacement jumps, a modified activation function with more physical meaning was introduced at the cohesive points (Figure \ref{fig:lr}). This adjustment allows for a delayed opening of cohesive elements.

Afterwards, the performance of the proposed PRNNs was evaluated. The three networks were trained with data from non-monotonic, non-proportional (GP-based) curves and tested on the same type of curves, and on proportional, non-monotonic curves with one cycle of unloading.

When tested on GP-based curves, the results showed that while all three networks followed the general trend of the curves, $\text{PRNN}_2$ and $\text{PRNN}_3$ performed with significantly higher accuracy than $\text{PRNN}_1$. Additionally, $\text{PRNN}_1$ failed to accurately capture the loss of stiffness due to damage evolution when predicting on the non-monotonic, proportional dataset. Specifically, $\text{PRNN}_1$ unloaded with the initial linear stiffness, which is due to the network's preference towards not utilizing the cohesive points effectively. This limitation can be explained by the network's layout. When the cohesive points are implemented in the material layer together with the bulk points, the stress output of the network is given by a linear combination of both stresses coming from the bulk models and tractions coming from cohesive zone models. This layout of the $\text{PRNN}_1$ does not resemble the physics of the full-order solution, where only the bulk points contribute to the stress homogenization.

On the other hand, the modified architectures with the damage input to the bulk points do not have this problem. Adjusting the local strain by the damage input allows for a modified tangent stiffness matrix that is able to capture the decrease in stiffness during unloading. This highlights the significance of designing the network's architecture with the knowledge of the underlying material behavior to achieve more accurate predictions.

When tested on non-monotonic, proportional curves with one cycle of unloading, $\text{PRNN}_3$ outperformed $\text{PRNN}_2$. It was observed that with a simple linear dense connection between the damage and bulk points, unloading and reloading occurred along different paths. While the RVE was unloading on the global scale, some bulk points in the network experienced further loading. This phenomenon occured because the weighted sum of damage made the input to the bulk point change sign, causing these points to undergo further loading instead of unloading. This caused further plastic strain development during macroscopic unloading and led to the different slopes during unloading and reloading. The issue with the different unloading/reloading path was mitigated when damage was used as an amplifier in $\text{PRNN}_3$. This method ensured that the fictitious bulk points follow the global trend of unloading/reloading.

Lastly, $\text{PRNN}_3$ was tested on non-monotonic, proportional curves with two cycles of unloading. The network provided accurate predictions in this case as well, demonstrating a progressively decreasing stiffness in successive unloading/reloading phases. This significant result highlights the network's capability to capture damage evolution over time, and further reinforces that the PRNN with modifications to its architecture is capable to represent microscale damage.




% \newpage
\bibliographystyle{elsarticle-num}  % Specify the bibliography style
\bibliography{report}





% \end{thebibliography}
\end{document}